\documentclass[a4paper,12pt]{article}

\usepackage{latexsym}
\usepackage{amsthm}
\usepackage{amsmath}
\usepackage{amsfonts, epsfig, amsmath, amssymb, color}
\usepackage[cmtip,arrow]{xy}
\usepackage{pb-diagram, pb-xy}
\usepackage{amssymb,epsfig}

\newfont{\sdbl}{msbm9}
\newfont{\dbl}{msbm10 at 12pt}
\theoremstyle{definition}
\newcommand{\cl}{{{\cal L}}}
\newcommand{\cq}{{{\cal Q}}}
\newcommand{\ch}{{{\cal H}}}
\newcommand{\cc}{{{\cal C}}}
\newcommand{\cn}{{{\cal N}}}
\newcommand{\ck}{{{\cal K}}}
\newcommand{\cj}{{{\cal J}}}
\newcommand{\cf}{{{\cal F}}}
\newcommand{\cm}{{{\cal M}}}
\newcommand{\co}{{{\cal O}}}
\newcommand{\ca}{{{\cal A}}}

\newcommand{\cg}{{{\cal G}}}
\newcommand{\xo}{{\mbox{\em \r{X}}}}
\newcommand{\eo}{{\mbox{\em \r{E}}}}
\newcommand{\dpp}{{\mbox{\dbl P}}}
\newcommand{\dz}{{\mbox{\dbl Z}}}

\newcommand{\dn}{{\mbox{\dbl N}}}

\newcommand{\sdz}{{\mbox{\sdbl Z}}}
\newcommand{\sdn}{{\mbox{\sdbl N}}}

\newcommand{\ord}{\mathop{\rm ord}\nolimits}

\newcommand\limproj{\mathop{\underleftarrow{\lim}}}
\newcommand\liminj{\mathop{\underrightarrow{\lim}}}

\makeatother
\newtheorem{defin}{Definition}
\newtheorem{notation}{Notation}
\newtheorem{ex}{Example}
\newtheorem{rem}{Remark}
\theoremstyle{plain}
\newtheorem{prop}{Proposition}

\newtheorem{theo}{Theorem}
\newtheorem{lemma}{Lemma}
\newtheorem{corol}{Corollary}

\newcounter{Figure}

\newcommand{\Proof}{{\noindent \bf Proof}}
\newcommand{\lto}{\longrightarrow}
\newcommand{\Spec}{\mathop{ \rm Spec}}
\newcommand{\Proj}{\mathop{ \rm Proj}}
\newcommand{\Ker}{\mathop{\rm Ker}}

\newcommand{\I}{\mathop{\rm Im}}

\title{Formal punctured ribbons and two-dimensional local fields}
\author{Herbert Kurke, \quad Denis Osipov, \quad Alexander Zheglov
\footnote{the second and the third authors are supported by RFBR
grant no. 08-01-00095, by INTAS grant 05-1000008-8118; besides the
second author is supported by grant of Leading Scientific Schools
no. 1987.2008.1, by grant of Russian Science Support Foundation,  by
a program of President of RF for supporting of young russian
scientists (grant no. MK-864.2008.1), and the third author is
supported by grant of Leading Scientific Schools no. 4578.2006.1,
and by grant of National Scientific Projects no. 2.1.1.7988. } }

\date{}

\begin{document}

\maketitle

\begin{abstract}
We investigate formal ribbons on curves. Roughly speaking, formal
ribbon is a family of locally linearly compact vector spaces on a
curve. We establish a one-to-one correspondence between formal
ribbons on curves plus some geometric data and some subspaces of
two-dimensional local field.
\end{abstract}

\section{Introduction}

The aim of this paper is to obtain an appropriate generalization of
the Krichever map for algebraic surfaces.

We recall that in the classical $1$-dimensional case it is a
one-to-one correspondence between integer projective curves over a
field $k$ plus line bundles (or torsion free sheaves if the curve is
singular) plus some additional data (a distinguished point $p$ of
the curve plus a formal local parameter at $p$, and a formal
trivialization at $p$ of the sheaf) and Schur pairs, i.e. pairs  of
$k$-subspaces $(W,A)$ of the  vector space $V=k ((z))$ satisfying a
Fredholm condition with respect to the subspace $V_+ =k[[z]]$ (i.e.
the complex $W\rightarrow V/V_+$ as well as the complex $A
\rightarrow V/V_+$ has to be Fredholm) such that $A$ is a
$k$-subalgebra of $V$ and $A \cdot W \subset W$.

Parshin and Osipov established the Krichever correspondence  in
higher dimensions (see \cite{Pa}, \cite{Os}, also \cite{Pa1},
\cite{ZhO}). In the 2-dimensional case it starts with a "flag" $(X
\supset C\supset p)$ (where
 $X$ is a  projective algebraic surface over a field $k$,
$C$ is an ample curve, $p$ is a $k$-point,  $X$ and $C$ are smooth
at $p$), a vector bundle $\cf$ of rank $r$ on $X$  plus formal
trivialization $e_p$ of $\cf$ at $p$, and formal local parameters
$u, t$ at $p$. By these data this correspondence associates the
$k$-subalgebra $A$ of $V=k ((u))((t))$ and $k$-subspace $W$ of
$V^{\oplus r}$ with Fredholm condition for all $i$ for $(A \cap
t^iV)/ (A \cap t^{i+1} V)$ as a subspace of $k((u))$ with respect to
$k[[u]]$, and with Fredholm condition for all $i$ for $(W \cap t^i
V^{\oplus r})/ (W \cap t^{i+1} V^{\oplus r})$ as a subspace of
$k((u))^{\oplus r}$ with respect to $k[[u]]^{\oplus r}$. $A$ is the
image of the structure sheaf $\co_X$, and $W$ is the image of the
sheaf $\cf$. If $X$ is Cohen-Macaulay surface, $C$ is an ample
Cartier divisor on $X$, then the pair $(A,W)$ contains all
information about $(X, C, p, \cf, e_p, u,t)$,
see~\cite[theorem~4]{Pa} and~\cite[theorem~6]{Os}.

However there was a problem, because contrary to the 1-dimensional
case it is not true that any such pair of subspaces comes from
geometric data. To solve this problem we introduce another type of
geometric objects which we call "ribbons" (or more exactly "formal
ribbons", but we will omit in the sequel the word "formal"). This
terminology comes from \cite{Ka}, where a similar object was
defined\footnote{More precisely, our ribbons are more general: the
ribbons from \cite{Ka} are $(C, \ca_0 )$ in our terminology}. We
decompose the Krichever map into the composition of the following
maps
$$
\left \{
\begin{array}{c}
\mbox{geometric data}\\
(X,C,p,\cf, e_p, u,t )
\end{array}
\right \}
\subset
\left \{
\begin{array}{c}
\mbox{geometric data}\\
\mbox{on ribbons}
\end{array}
\right \}
\mapsto
\left \{
\begin{array}{c}
\mbox{pairs of subspaces } (W,A)\\
\mbox{ with Fredholm conditions} \\
\end{array}
\right \}
$$
Ribbons are ringed spaces which are, on the one hand side,  more
general as the notion of "formal scheme" of Grothendieck, on the
other hand side, they have some extra features. We explain them
exactly in section~\ref{rib}.

In section~\ref{coh} we clarify the cohomology of sheaves, which we
call ind-pro-quasicoherent sheaves on a ribbon. We investigate the
coherence property of ribbons.

In section~\ref{Pic} we clarify the structure of the Picard group of
a ribbon.

In section~\ref{KPmap}  we establish a one-to-one correspondence
between the classes of isomorphic "geometric data" (punctured ribbon
plus torsion free sheaf on it plus some extra data) and the "Schur
pairs" $(A,W) \subset (V, V^{\oplus r})$, were $A$ is a
$k$-subalgebra of $V$, and $A \cdot W \subset W$, satisfying
Fredholm conditions for the subquotients (as explained above).

We computed also several examples to illustrate the general theory.

We note that  families of Tate  spaces (i.e. of locally linearly
compact vector spaces) were studied also in~\cite{Dr}.

We think that the ribbons and geometric data on them, which are
introduced in this paper, will help to find geometric solutions of
generalizations of Parshin's two-dimensional KP-hierarchy,
see~\cite{Pa1}, \cite{Zh}.

\bigskip

{\bf Aknowledgements.} This research was done at the
Mathematisches Forschungsinstitut Oberwolfach during a stay within the
Research in Pairs Programme from January 28 - February 10, 2007.
We
would like to thank the MFO at Oberwolfach for the excellent working
conditions.

We are grateful to A.N. Parshin for his interest to this work.

\section{The category of ribbons}
\label{rib}
\subsection{Definition of a ribbon.}

Let $S$ be a Noetherian base scheme.
\begin{defin} \label{ribbon}
A ribbon $(C, \ca)$ over $S$ is given by the following data.
\begin{enumerate}
\item \label{a}
 A flat family of reduced algebraic curves $\tau :  C\rightarrow S$.
\item \label{b} A sheaf $\ca$ of commutative $\tau^{-1}\co_S$-algebras
on $C$.
\item \label{c} A descending sheaf filtration $(\ca_i)_{i\in \sdz}$ of $\ca$ by
$\tau^{-1}\co_S$-submodules which satisfies  the following axioms:
\begin{enumerate}
\item  \label{i} $\ca_i\ca_j \subset \ca_{i+j}$, $1\in \ca_0$ (thus $\ca_0$
is a subring, and for any $i \in \dz$ the sheaf $\ca_i$ is a
$\ca_0$-submodule);
\item \label{ii}
$\ca_0/\ca_1$ is the structure sheaf $\co_C$ of $C$;
\item
\label{iii} for each $i$ the sheaf $\ca_i/\ca_{i+1}$ (which is a
$\ca_0/\ca_1$-module by (\ref{i})) is a coherent sheaf on $C$, flat
over $S$, and for any $s \in S$ the sheaf $\ca_i/\ca_{i+1}
\mid_{C_s}$ has no coherent
subsheaf with finite support, and is isomorphic to $\co_{C_S}$ on a dense open set;
\item
\label{iv}  $\ca =\liminj\limits_{i \in \sdz} \ca_i$, and
$\ca_i=\limproj\limits_{j>0}\ca_i/\ca_{i+j}$ for each $i$.
\end{enumerate}
\end{enumerate}
\end{defin}

\begin{rem}
It follows from  (\ref{iii}) of the definition that if $C_s$ (for $s
\in S$) is an irreducible curve, then the sheaf $\ca_i /\ca_{i+1}
\mid_{C_s}$ is a torsion free sheaf on $C_s$ for any $i \in\dz$.
\end{rem}

\begin{notation}
For simplicity we denote a ribbon $(C, \ca)$ over $\Spec R$, where
$R$ is a ring, as a ribbon over $R$.
\end{notation}

\begin{ex}
\label{ex1} If $X$ is an algebraic surface over a field $k$, and
$C\subset X$ is a reduced effective Cartier divisor, we obtain a
ribbon $(C, \ca )$ over $k$, where
$$ \ca := \co_{\hat{X}_C}(*C)=\liminj\limits_{i \in \sdz} \co_{\hat{X}_C}(-iC)
= \liminj\limits_{i \in \sdz} \: \limproj\limits_{j \ge 0} J^{i}/
J^{i+j}
$$
$$
 \ca_i := \co_{\hat{X}_C}(-iC) = \limproj\limits_{j \ge 0} J^{i}/
J^{i+j}  \mbox{,} \quad i \in \dz  \mbox{,}$$
 where $\hat{X}_C$ is the formal scheme which is the completion of $X$ at
 $C$,
 and $J$ is the ideal sheaf of $C$ on $X$ (the sheaf $J$ is an invertible
 sheaf).
\end{ex}

\begin{prop}
\label{prop1}
\begin{enumerate}
\item
\label{1)} For any $i \ge 0$ the ringed space
$X_i=(C,\ca_0/\ca_{i+1})$, $(i\ge 0)$ is a scheme, which is flat over $S$.
\item
\label{2)} For any $j \in \dz$ and any $i \ge 0$ the sheaf
$\ca_j/\ca_{j+i+1}$ is a coherent sheaf on $X_i$, which is a flat
sheaf over $S$.
\item
\label{3)} If $X_{\infty}=(C,\ca_0)$, then $X_{\infty}$ is a locally
ringed space, and  we have closed embeddings of schemes
$$
X_0 \subset X_1 \subset X_2 \subset \ldots  X_i\subset X_{i+1}
\subset \ldots
$$
such that $X_{\infty}=\liminj\limits_{i \ge 0} \: X_i$ in the
category of locally ringed spaces.
\end{enumerate}
\end{prop}
\Proof. We prove the first statement of the proposition.

At first, we show that $X_i$ are locally ringed spaces. By
definition, we have that $X_0$ is the scheme $(C,\co_C)$. Therefore
$X_0$ is a locally ringed space. We have that for every $i \ge 0$
the subsheaf $\ca_i/\ca_{i+1} \subset \ca_0/\ca_{i+1}$ is a
nilpotent ideal sheaf because of~(\ref{i}) of
definition~\ref{ribbon}. We consider the following exact triple of
sheaves on $C$:
\begin{equation}
\label{*} 0 \lto \ca_i/\ca_{i+1} \lto \ca_0/\ca_{i+1}
\arrow{e,t}{\pi_i} \ca_0/\ca_{i}   \lto 0 \mbox{.}
\end{equation}
For each point $P \in C$ we consider the stalks at $P$ of each sheaf
from this sequence. We obtain the following exact sequence:
\begin{equation}
\label{2*} 0 \arrow{e} (\ca_i/\ca_{i+1})_P \arrow{e}
(\ca_0/\ca_{i+1})_P  \arrow{e,t}{(\pi_{i})_P} (\ca_0/\ca_{i})_P
\arrow{e} 0 \mbox{.}
\end{equation}
 We  apply now induction arguments on $i$. By induction hypothesis, we assume
 that the ring
$(\ca_0/\ca_{i})_P$ is a local ring with the maximal ideal $\cal M$.
Let $ \cm'$ be the ideal $\pi_{i}^{-1}(\cm )$. Then this ideal  is a
unique maximal ideal in $(\ca_0/\ca_{i+1})_P$. Therefore this ring
is a local ring. Indeed, if $a\in (\ca_0/\ca_{i+1})_P \backslash
\cm'$, then $a$ must be invertible in the ring
$(\ca_0/\ca_{i+1})_P$, since $(\pi_i)_P (a)$ is invertible in the
ring $(\ca_0/\ca_{i})_P$, and $(\ca_i/\ca_{i+1})_P$ is a nilpotent
ideal in the ring $(\ca_0/\ca_{i+1})_P$.

Second, we show that there are natural morphisms $X_i
\arrow{e,t}{\tau_i} S$ of locally ringed spaces for each $i \ge 0$,
and  that these morphisms are flat. We apply induction on $i \ge 0$.
For every $i \ge 0$ the morphism $\tau_i$ consists of the
topological morphism $\tau : C\rightarrow S$
 and of morphism of sheaves
 $$\tau_i^{\sharp} \; : \; \co_S\rightarrow
\tau_*(\ca_0/\ca_{i+1}), \qquad \mbox{where} $$
$$\tau_i^{\sharp}(U) \; : \; \co_S(U)\ni a \longmapsto a \cdot 1 \in
\ca_0/\ca_{i+1}(\tau^{-1}(U))$$ for each open subset $U \subset S$.
For each $P \in C$ the morphism
$$(\tau_i^{\sharp})_P \; : \; (\co_{S})_{\tau(P)} \lto
(\ca_0/\ca_{i+1})_P$$ is a local morphism, because its composition
with the morphism $(\pi_ i)_P$ is a local morphism by induction
hypothesis.

Now for every $i \ge 0$ the morphism $\tau_i$ is a flat morphism by
standard results on flat modules (see e.g. \cite[ch. 2, \S 3]{Ma}),
because for each $P \in C$ the $(\co_{S})_{\tau(P)}$-modules
$(\ca_i/\ca_{i+1})_P$ and $(\ca_0/\ca_{i})_P$  are flat
$(\co_{S})_{\tau(P)}$-modules by induction hypothesis on $i$ and
by~(\ref{iii}) of definition~\ref{ribbon}. Therefore we obtain from
exact sequence~(\ref{2*}) that $(\ca_0/\ca_{i+1})_P$ is a flat
$(\co_{S})_{\tau(P)}$-module.

At third, we show that a locally ringed space $X_i$ is scheme for
each $i \ge 0$.  We consider any affine open subset $U \subset C$.
The sequence~(\ref{*}) leads to the following exact triple:
\begin{equation} \label{eq}
0\rightarrow \ca_i/\ca_{i+1}(U) \longrightarrow \ca_0/\ca_{i+1}(U)
\overset{\pi}{\longrightarrow} \ca_0/\ca_{i}(U)\longrightarrow 0.
\end{equation}
This sequence is an exact sequence, because the sheaf
$\ca_i/\ca_{i+1}$ is a coherent sheaf on $C$, and $U$ is an affine
set. We have that $\ca_0/\ca_{i+1}(U)$ and $\ca_0/\ca_{i}(U)$ are
rings, and we are going to show that $(U,(\ca_0/\ca_{i+1})|_U)\simeq
\Spec (\ca_0/\ca_{i+1}(U))$.

It is clear that the topological space
$ \underline{\Spec}(\ca_0/\ca_{i+1}(U))=U  \mbox{.}$
Using that $\ca_i/\ca_{i+1}(U)$ is a nilpotent ideal in the ring
$\ca_0/\ca_{i+1}(U)$, from exact sequence~(\ref{eq}), by induction
on $i$ we obtain that the identical map on the ring
$\ca_0/\ca_{i+1}(U)$ induces a well-defined morphism of sheaves on
$U$ :
$$
\gamma \: : \: \widetilde{\ca_0/\ca_{i+1}(U)} \longrightarrow
(\ca_0/\ca_{i+1}) \mid_U \mbox{,}
$$
where for any $\ca_0/\ca_{i+1}(U)$-module $ N$ by $\widetilde{N}$ we
denote the corresponding quasicoherent sheaf on
$\Spec(\ca_0/\ca_{i+1}(U))$. The map $\gamma$ is an isomorphism,
since it follows from the following exact diagram of sheaves on $U$:
$$
\begin{array}{ccccccc}
0\longrightarrow &\widetilde{\ca_i/\ca_{i+1}(U)}&\longrightarrow
&\widetilde{\ca_0/\ca_{i+1}(U)}& {\rightarrow} &\widetilde{\ca_0/\ca_{i}(U)}&
\longrightarrow 0\\
&\downarrow && \downarrow &&\downarrow & \\
0\longrightarrow & (\ca_i/\ca_{i+1})|_U & \longrightarrow &
(\ca_0/\ca_{i+1})|_U & \rightarrow & (\ca_0/\ca_{i+1})|_U &
\longrightarrow 0
\end{array}
$$
The left vertical arrow in this diagram  is an isomorphism
by~(\ref{iii}) of definition~\ref{ribbon}. The right vertical arrow
is an isomorphism by induction on $i$. Therefore, the middle
vertical arrow $\gamma$ is also an isomorphism. Thus we proved that
$X_i$ is a scheme for each $i \ge 0$. It finishes the proof of the
first statement of the proposition.

We prove now the second statement of the proposition. As above, the
proof is by induction on $i$. We have the following exact triple of
$\co_{X_i}$-modules:
$$
0 \longrightarrow \ca_{j+i}/\ca_{j+i+1} \longrightarrow
\ca_j/\ca_{j+i+1} {\longrightarrow} \ca_j/\ca_{j+i} \longrightarrow
0.
$$

By definition, the sheaf $\ca_{j+i}/\ca_{j+i+1}$ is a coherent sheaf
on $X_i$, and a flat sheaf over $S$.  The sheaf $\ca_j/\ca_{j+i}$ is
a coherent $\co_{X_{i-1}}$-module sheaf, and  a flat sheaf over $S$
by the induction hypothesis. Therefore, this sheaf is  also a
coherent $\co_{X_i}$-module sheaf, because both module structures
coincide on this sheaf. Thus, the sheaf $\ca_j/\ca_{j+i+1}$ is a
coherent sheaf by \cite[prop. 5.7]{Ha} and flat over $S$ by
\cite[prop. 9.1]{Ha}. We proved the second statement of the
proposition.

The third statement of the proposition easily follows from exact
sequence~(\ref{*}).
\begin{flushright}
$\square$
\end{flushright}

\begin{defin} \label{defisom}
\begin{enumerate}
\item
A  morphism $\varphi$ of ribbons over $S$
$$\varphi :(C,\ca
)\rightarrow (C',\ca')$$ is a morphism of ringed spaces over $S$
that preserves the filtrations, i.e. we have for the map
$\varphi^{\sharp} : \ca' \to \varphi_*(\ca)$, for any $i \in \dz$
$$ \varphi^{\sharp}(\ca'_i)\subset
\varphi_*(\ca_i) \mbox{.}$$
\item
An  isomorphism of ribbons is a morphism that has right and left
inverse.
\end{enumerate}
\end{defin}

\subsection{Base change.}

\begin{notation} We will also denote the ribbon $(C,\ca )$ by
$\xo_{\infty}$.
\end{notation}

For a ribbon $\xo_{\infty} = (C, \ca)$ over $S$, and a morphism
$\alpha: S' \longrightarrow S$ of Noetherian schemes we define a
base change ribbon $\xo_{\infty}' = (C', \ca')$ over $S'$
in the following way:
$$
C' :=C \times_S S'  \mbox{,}
$$
$$
\ca'_j :=\limproj\limits_{i\ge 1}
(\ca_j/\ca_{j+i})\boxtimes_{\co_S}\co_{S'}
$$
for any $j \in \dz$. From statement~\ref{2)} of
proposition~\ref{prop1} we have for any $j \in \dz$, any $i \ge 0$
$$(\ca_{j+1}/\ca_{j+i+1}) \boxtimes_{\co_S}\co_{S'}\subseteq
(\ca_{j}/\ca_{j+i+1}) \boxtimes_{\co_S}\co_{S'} \mbox{.}$$ Therefore
we have
$$
\ldots \subset  \ca'_{j+1} \subset\ca'_j \subset \ca'_{j-1} \subset
\ldots
$$
and we define
$$
\ca' :=\liminj_{i \in \sdz}\ca'_i.
$$
By definition we have $\ca'_j/\ca'_{j+1}=(\ca_j/\ca_{j+1})
\boxtimes_{\co_S}\co_{S'}$ for any $j \in \dz$, and all axioms from
definition of ribbon are satisfied.

\begin{prop}
\label{prop2} For the base change $\alpha : S' \longrightarrow S$ we
have that for the base change ribbon $\xo'_{\infty} = (C', \ca')$
the following properties are satisfied.
\begin{enumerate}
\item $X'_i=X_i\times_SS'$ for any $i \ge 0$
\item
$\ca'_j/\ca'_{j+i+1}=(\ca_j/\ca_{j+i+1}) \boxtimes_{\co_S} \co_{S'}$
for any $j \in \dz$ and any $i \ge 0$
\end{enumerate}
\end{prop}
\Proof. The proof is clear from definition of a ribbon and
proposition~\ref{prop1}.
\begin{flushright}
$\square$
\end{flushright}

\section{Coherent sheaves on a ribbon.} \label{coh}
\subsection{Ind-pro-quasicoherent sheaves.}

\begin{defin} \label{ipqs}
Let $\xo_{\infty} = (C, \ca )$ be a ribbon, and $\cf$ a sheaf of
$\ca$-modules. We will call $\cf$ {\it ind-pro-coherent
(ind-pro-quasicoherent)} on $\xo_{\infty}$ if it has a descending
sheaf filtration $(\cf_i)_{i\in \sdz}$  with the following
properties.
\begin{enumerate}
\item
\label{s1} $\ca_i\cf_j\subseteq\cf_{i+j}$.
\item
\label{s2} $\cf_j/\cf_{j+1}$ is a coherent (quasicoherent)
$\co_C$-module for all $j$.
\item
\label{s3} $\cf_i=\limproj\limits_j\cf_i/\cf_{i+j}$.
\item
\label{s4}
$\cf =\liminj\limits_i\cf_i$.
\end{enumerate}
\end{defin}

We recall that projective system  $(D_i, i \in \dn)$ of abelian
groups with transition maps $\phi_{i'i}$ satisfies the {\em
ML-condition} (the Mittag-Lefler condition), iff for every $i \in
\dn$ the descending family of subgroups $\{ \phi_{i'i} (D_{i'})
\subset D_i \mid i' \ge i \in \dn \}$ will stabilize.

We will need the following lemma, which is easy to prove,
using~\cite[prop.~9.1]{Ha}.
\begin{lemma} \label{lemma1}
If
$$
0 \longrightarrow (K_i) \longrightarrow (A_i) \longrightarrow (B_i)
\longrightarrow (C_i) \longrightarrow 0
$$
is an exact sequence of projective systems of abelian groups with
respect to $\dn$, and projective systems $(K_i, i\in \dn)$ and
$(A_i, i \in \dn)$ satisfy the ML-condition, then the induced
sequence of projective limits
$$
0 \longrightarrow  \limproj\limits_{i \in \sdn} K_i \longrightarrow
\limproj\limits_{i \in \sdn} A_i \longrightarrow
\limproj\limits_{i\in \sdn} B_i  \longrightarrow \limproj\limits_{i
\in \sdn} C_i \longrightarrow 0
$$
is also exact.
\end{lemma}

\begin{prop}
\label{prop3}
Let $\xo_{\infty}=(C,\ca )$ be a ribbon and $\cf$ an
ind-pro-quasicoherent sheaf on $\xo_{\infty}$. Then  we have the
following.
\begin{enumerate}
\item
\label{i1}
 $\cf_i/\cf_{i+j+1}$ is a quasicoherent $\co_{X_j}$-module for any $j \ge 0$, $i \in \dz$.
\item
\label{i2} We have   that $\cf_i(U)/\cf_j(U)\rightarrow
(\cf_i/\cf_j)(U)$ is an isomorphism for all $i<j$ and for any affine
open subset $U\subset C$.
\item
\label{i3} If $\xo_{\infty}$ is a ribbon over  an  Artinian ring,
then  for any affine open subset $U\subset C$ we have
$H^1(U,\cf_i)=H^1(U,\cf )=0$.
\end{enumerate}
\end{prop}
\Proof. The proof of statement~\ref{i1} of the proposition is
analogous to the proof of statement \ref{2)} of proposition
\ref{prop1}.

We prove statement~\ref{i2} of the proposition. We always have an
exact sequence
$$
0\rightarrow \cf_j(U)\rightarrow \cf_i(U) \rightarrow
(\cf_i/\cf_j)(U) \mbox{,}
$$
and we have  exact sequences for $i<j<k$
$$
0\rightarrow (\cf_j/\cf_k)(U)\rightarrow (\cf_i/\cf_k)(U)\rightarrow
(\cf_i/\cf_j)(U)\rightarrow 0 \mbox{,}
$$
since, by statement~\ref{i1}, $\cf_j/\cf_k$ is a quasicoherent sheaf
of $\co_{X_{k-j-1}}$-modules.

Now since $\cf_i(U)=\limproj\limits_{k \ge i}(\cf_i/\cf_k)(U)$ and
all maps $(\cf_j/\cf_{k+1})(U)\rightarrow (\cf_j/\cf_k)(U)$ are
surjective, we also have surjections $\cf_i(U)\rightarrow
(\cf_i/\cf_j)(U)$ (see lemma~\ref{lemma1}).

We prove statement~\ref{i3} of the proposition. Since $C$ is a curve
over an Artinian ring, every open subset of an affine open set $U$
is again affine. We take an embedding $\cf_i\hookrightarrow W$ into
a flabby sheaf, then $H^1(U,\cf_i)$ is the cokernel of
$W(U)\rightarrow (W/\cf ) (U)$, and we have to show that any section
of $(W/\cf )(U)$ lifts to a section of $W(U)$.

Since the underlying space $U$ is Noetherian, we have a largest open
set $U'\subseteq U$  where a lifting $w'$ of the given section
exists. We will show that the assumption $U'\subsetneqq U$ leads to
a contradiction. Assume $p\in U\backslash U'$, then we find a
neighbourhood $U''\subset U$ of $p$ and a lifting $w''$ on $U''$ of
the given section. If $U'\cap U''=\varnothing$ we could extend
$(U',w')$ to $(U'\cup U'', w' \mbox{ on } U', w'' \mbox{ on } U'')$.
If $U'\cap U''\ne \varnothing$, we get a section $a=w' -w''$ of
$\cf_i(U'\cap U'')$.

We {\em claim} that $\cf_i(U')\oplus \cf_i(U'')\rightarrow
\cf_i(U'\cap U'')$ is surjective, so we can write $a=a'-a''$ with
$a'\in \cf_i(U')$, $a''\in \cf_i(U'')$. Then $w|_{U'}=w'-a'$ and
$w|_{U''}=w''-a''$ would give a lifting to $U'\cup U''$, hence $U'$
was non maximal.

{\em Proof of the claim.} We have an exact sequence of projective
systems
$$
\begin{array}{ccccccc}
&\downarrow &&\downarrow &&\downarrow & \\
0\rightarrow &{\cf_i/\cf_{j+1}(U'\cup U'')}&\rightarrow &{\cf_i/\cf_{j+1}(U')
\oplus \cf_i/\cf_{j+1}(U'')}& {\rightarrow} &{\cf_i/\cf_{j+1}(U'\cap U'')}&\rightarrow 0\\
&\downarrow &&\downarrow &&\downarrow & \\
0\rightarrow & (\cf_i/\cf_{j})(U'\cup U'') & \rightarrow & (\cf_i/\cf_{j})(U')
\oplus (\cf_i/\cf_{j})(U'') & \rightarrow & (\cf_i/\cf_{j})(U'\cap U'') & \rightarrow 0 \\
&\downarrow &&\downarrow &&\downarrow &
\end{array}
$$
where all transition maps are surjective. Thus the projective limit
stays exact (see lemma~\ref{lemma1}). For $\cf$ the assertion
follows since cohomology commute with $\liminj$.
\begin{flushright}
$\square$
\end{flushright}

\begin{corol} \label{col1}
Let $\xo_{\infty} = (C, \ca)$ be a ribbon over  $A$, where $A$ is an
Artinian  ring. Let $\cf$ be an ind-pro-quasicoherent sheaf on
$\xo_{\infty}$, and $C$ be a projective curve over $A$.
\begin{enumerate}
\item
If $C = U_1 \cup U_2$, where  $U_1$ and $U_2$ are affine open
subsets, then we have an exact sequence
$$
0 \rightarrow H^0(C,\cf ) \rightarrow H^0(U_1,\cf )\oplus
H^0(U_2,\cf ) \rightarrow H^0(U_1 \cap U_2, \cf ) \rightarrow
H^1(C,\cf )\rightarrow 0 \mbox{.}
$$
\item
If $\cf$ is an ind-pro-coherent sheaf, then
$$H^*(C,\cf )=\liminj\limits_i\limproj\limits_{j\ge i}H^*(X_{j-i}, \cf_i/\cf_{j+1}) \mbox{.}$$
\end{enumerate}
\end{corol}
\Proof. The first assertion of this corollary is the Mayer-Vietoris
exact sequence, due to assertion~\ref{i3} of proposition
\ref{prop3}, because $U_1$ and $U_2$ are affine sets.

 We prove now  the second
assertion of this corollary. We note that for any $j \in \dz$ a
projective system ($H^0(C, \cf_j / \cf_{j+i}), i \in \dn$) satisfies
the ML-condition, because $H^0 (C, \cf_j / \cf_{j+i})$ is an
Artinian $A$-module for any $i, j$, and the maps in projective
system are the maps of $A$-modules.

We note that, since $C$ is a curve over an Artinian ring, there are
some affine open subsets $U_1$ and $U_2$ of $C$ such that $C= U_1
\cup U_2$. For any fixed $j \in \dz$ a projective system
$(H^0(U_1,\cf_j / \cf_{j+i}) \oplus H^0(U_2, \cf_j / \cf_{j+i} ), i
\in \dn)$ satisfies the ML-condition because of assertion~\ref{i2}
of proposition~\ref{prop3}.

Now, since the cohomology commutes with direct limits,  the second
assertion of this corollary follows from the first one, using
lemma~\ref{lemma1}.
\begin{flushright}
$\square$
\end{flushright}

\subsection{Coherence property}

\begin{rem} \label{rem1}
The sheaf $\ca$ may be not coherent in the usual sense (due to H.
Cartan, see \cite{Se}).

We recall that a sheaf $\cf$ of $\ca$-modules on a topological space
$X$ is coherent if it satisfies the following two properties.
\begin{enumerate}
\item $\cf$ is locally of finite type, i.e. for any point $x\in X$
there exist an open $U\ni x$ and finite number of sections $s_1,
\ldots ,s_p\in \cf (U)$ such that for any $y\in U$ the stalk $\cf_y$
is generated by the images of $s_1, \ldots ,s_p$ over $\ca_y$.
\item
The sheaf $\ck=ker ((\ca |_U)^{\oplus q}\overset{(f_1,\ldots
,f_q)}{\longrightarrow} (\cf |_U))$, where $f_i\in \cf (U)$ for an
open $U$, is locally of finite type. Here the map
$\overset{(f_1,\ldots ,f_q)}{\longrightarrow}$ maps an element
$(a_1,\ldots ,a_q)$ to $\sum a_if_i$.
\end{enumerate}

The sheaf $\ca$ is called coherent if it is coherent as $\ca$-module.
\end{rem}

 Let's consider the following
ringed space: $(C, \co_C ((t))^{Q} )$, where $C$ is a reduced
algebraic curve over a field $k$, $Q \in C$ is a closed point, and
 the sheaf $\co_C ((t))^{Q}$ is defined by
$$
\co_C ((t))^{Q}(U):= \{ \sum_{i =l}^{\infty}c_it^i, \mbox{ where }
c_i\in \co_C(U) \mbox{ for $i\ge 0$ } \mbox{ and } c_i\in \cj_Q(U)
\mbox{ for $i<0$ }\} \mbox{,}
$$
where $\cj_Q$ is the ideal sheaf of the point $Q$. Clearly, this is
a sheaf, and $(C, \co_C ((t))^{Q})$ is a ribbon over the field $k$.
This sheaf is an analogue of the sheaf $\co_X((t))^{\sqrt{}}$ from
\cite{KV1}.

\begin{ex}
\label{ex2}
 This is an example of non-coherent sheaf $\ca$ of a ribbon.

 Let $C$ be a plane affine singular cubic curve given
 by the equation $y^2=x^2(x+1)$ over a field $k$, $Q\in C$ is a closed point $x=y=0$.
We show that the sheaf  $\ca=\co_C ((t))^{Q}$ is not coherent.

If it were coherent, then, by definition, for each $q \ge 1$ and
$f_1,\ldots f_q\in \ca (U)$ the sheaf $\ck=ker ((\ca |_U)^{\oplus
q}\overset{(f_1,\ldots ,f_q)}{\longrightarrow} (\ca |_U))$  must be
locally of finite type.  We  take $U\ni Q$, $q=2$, and let $f_1,
f_2$ be the images of $x,y$ in $\co_C(U)$. Let $V\subset U$, $Q\in
V$ be an open set such that $\ck (V)$ is finitely generated in each
point.

We consider an element $(b_1,b_2)  \in \ck (V)$ such that $b_1, b_2$
are the images of $-y, x$ in $\co_C(U)$. Then $(b_1, b_2) \in
\cj_Q((t))^{\oplus 2}(V)$, but $(b_1, b_2)\notin
\cj_Q^{2}((t))^{\oplus 2}(V)$. We note that elements
$(b_1t^m,b_2t^m) \in \ck (V)$ for any  $m \in \dz$, and also satisfy
this condition.

We note that for each $(a_1,a_2)\in \ck (V)$ we have $a_i=\sum
a_{ij}t^j$, where $a_{ij}\in \cj_Q (V)$. Indeed, we must have
$f_1a_{1j}+f_2a_{2j}=0$ for all $j$, and this equality holds only if
$a_{ij}$ are polynomials in $f_1,f_2$ without free terms, i.e.
belong to the ideal $\cj_Q(V)$.

Let $g_1,\ldots g_l$ be generators of $\ck (V)$. Let they have
orders $(q_i,q'_i)$, $i=1,\ldots l$, where the order of an element
$a\in \ca (V)$ is  equal to the degree (with respect to $t$) of the
lowest term of $a$. For each $m \in \dz$  we must have
\begin{equation} \label{eq1}
(b_1t^m,b_2t^m)=\sum_{i=1}^l w_{im}g_i
\end{equation}
 with $w_{im}\in \ca (V)  $. If $M = \min \{ q_1, \ldots, q_l,  q_1' \ldots, q_l'\}$,
 then all coefficients of
$t^j$ with $j<M$ on the right hand side of~formula (\ref{eq1}) must
belong to $\cj_Q^{2}(V)^{\oplus 2}$ for each $m$. But if $m \ll 0$
then there will be coefficients of $t^j$ with $j<M$ on the left hand
side of formula~(\ref{eq1}) that do not belong to
$\cj_Q^{2}(V)^{\oplus 2}$ (and the same is true for their images in
the stalk of $Q$). We have a contradiction.

The same arguments show that the ideal $\cj_Q(V)((t))\subset \ca
(V)$ is not finitely  generated, i.e. the ring $\ca (V)$ is not
Noetherian.
\end{ex}

For convenience, we introduce also the following definition.

\begin{defin}
The sheaf of rings $\cf$ on a topological space $X$ is called weakly
Noetherian, if there exists an open affine cover
$\{U_{\alpha}\}_{\alpha\in I}$ such that $\cf (U_{\alpha})$ is a
Noetherian ring for any $\alpha\in I$.
\end{defin}

\begin{ex}
\label{ex2a} This is an example of coherent, but not weakly
Noetherian sheaf $\ca$ of a ribbon.

We  consider the ringed space
 $(C, \ca =\co_C ((t))^{Q} )$, where $C$ is a reduced algebraic curve over a field $k$, $Q \in C$ is a smooth point.
 We will prove that the sheaf $\ca$ is a coherent sheaf of rings.
 To  prove that the sheaf $\ca$ is a coherent sheaf of rings, it is enough to prove that the sheaf
 $\ck$ from definition of coherence (see remark \ref{rem1} above) is locally of finite type.

We consider an open $U \subset C$. If $U \not\ni Q$, then we have
$(\ca |_U)^{\oplus q} \simeq (\co_C ((t))|_U)^{\oplus q}$ and
therefore for any affine open subset $V\subset U$ the ring $(\ca
|_U)(V)$ is Noetherian. Clearly, $\ck (V)=(\ck' (V))_t$ and $\ca
(V)=(\ca' (V))_t$, where
$$
\ck'=ker ((\ca' |_U)^{\oplus q}
\overset{(f_1t^{k},\ldots
,f_qt^{k})}{\longrightarrow} (\ca' |_U)), \mbox{\quad}
\ca'=\co_C[[t]]
$$
for sufficiently large $k$ (note that the definition of the sheaf
$\ck$ does not depend on changes $(f_1, \ldots ,f_q)\mapsto
(f_1t^{k},\ldots ,f_qt^{k})$). The locally ringed space $(C, \ca')$
is a Noetherian formal scheme (so, $\ca'$ is a coherent sheaf, see
\cite[ch.I, \S 10.10]{EGAI}), therefore $\ck'$ is locally of finite
type, i.e. for each point $P\in U$ there exists open $V\subset U$,
$P\in V$ and generators $(\beta_1, \ldots ,\beta_n)$ of $\ck'(V)$
over $\ca'(V)$ such that their images generate stalks $\ck'_x$ for
each $x\in V$. Clearly, $(\beta_1, \ldots ,\beta_n)$ are also
generators of the $\ca (V)$-module $\ck (V)$, and they also generate
stalks $\ck_x$ over $\ca_x$ for each $v \in V$, i.e. $\ck$ is
locally of finite type.

Let now $U \ni Q$, $f_1,\ldots , f_q\in \ca (U)$. Our sheaf $\ca$ is
a subsheaf of the sheaf \linebreak $\tilde{\ca}=\co_C ((t))$, and
the last sheaf is coherent, as we have proved above. We define the
sheaf
$$
\tilde{\ck}=ker ((\tilde{\ca} |_U)^{\oplus q}
\overset{(f_1,\ldots
,f_q)}{\longrightarrow} (\tilde{\ca} |_U)).
$$
It is locally of finite type, and $\ck$ is the subsheaf of $\tilde{\ck}$ as a sheaf of abelian groups.

For a given element $a = \sum\limits_j a_jt^j \in \tilde{\ca}(V)$,
$Q\in V$ we define its $Q$-order as follows:
$$
\ord_Q(a) = \left\{
\begin{array}{l}
\min \; \{j \; : \;   a_j\notin \cj_Q \}  \\
\infty \mbox{,} \quad \mbox{if for any} \quad j \quad  a_j \in \cj_Q
\mbox{.}
\end{array}
\right.
$$

Clearly, for any $a,b\in \tilde{\ca}(V)$ we have
$$\ord_Q(ab)=\ord_Q(a)+\ord_Q(b) \mbox{.}$$

For a given element $\alpha \in \tilde{\ca}^{\oplus q}(V)$, $Q\in V$
we define its $Q$-order  as a minimum of $Q$-orders of components of
$\alpha$, i.e.
$$
\ord_Q(\alpha) = \min \{ a_1, \ldots, a_q\} \quad \mbox{for} \quad
\alpha=(a_1, \ldots, a_q) \mbox{.}
$$

Let $\alpha_1, \ldots, \alpha_k$ be generators of the
$\tilde{\ca}(V)$-module $\tilde{\ck}(V)$, $V\ni Q$, such that their
images generate stalks $\tilde{\ck}_x$ for each $x\in V$. Without
loss of generality we can assume that $V$ is an affine open set,
such that the maximal ideal of the point $Q$ in $\co_C(V)$ is a
principal ideal $(y)$, $y\in \co_C(V)$. We can also assume that
$\alpha_1, \ldots, \alpha_k\in \ck (V)$, since otherwise we can
replace them by $\alpha_1t^{l_1}, \ldots \alpha_kt^{l_k}$. Since the
maximal ideal of the point $Q$ in $\co_C(V)$ is a principal ideal,
 we have $\alpha_i=y^{k_i}\alpha'_i$,
where $k_i\ge 0$ and $\ord_Q(\alpha'_i)<\infty$. We assume
$\alpha'_i\in \ck (V)$ again after multiplication them by some
powers of $t$. Obviously, the elements $\alpha'_i\in \ck (V)$ are
also generators of $\tilde{\ck}(V)$ and of stalks $\tilde{\ck}_x$
for each $x \in V$. So, we can assume that $\ord_Q(\alpha_i)=0$ for
any $1 \le i \le k$.

Without loss of generality we can assume that the first component of
$\alpha_1$ be of zero $Q$-order. Since the ring $\co_C(V)$ has
dimension $1$, we can change $\alpha_1, \ldots, \alpha_k$ by
$\alpha_1, \alpha_2+x_2\alpha_1, \ldots , \alpha_k+x_k\alpha_1$ for
some $x_2, \ldots ,x_k\in \ca(V)$
 such that the first components of elements $\alpha_2+x_2\alpha_1, \ldots , \alpha_k+x_k\alpha_1$ has
 infinite $Q$-order.
 If the $Q$-order of an element $\alpha_i+x_i\alpha_1$ is finite,
 we can again assume that it is zero, after multiplication him by an appropriate power of $t$.

The elements $\alpha_1, \alpha_2+x_2\alpha_1, \ldots ,
\alpha_k+x_k\alpha_1$ are again generators of $\tilde{\ck}(V)$ (and
of $\tilde{\ck}_x$ for each $x \in V$). They form a $k\times q$
matrix, whose entries lie in $\ca (V)$ (the $i$-th row is the
element $\alpha_i+x_i\alpha_1$). The corresponding $k\times q$
matrix of its $Q$-orders looks like
$$
\left (
\begin{array}{cccc}
0& *& \ldots &*\\
\infty &*& \ldots &*\\
\vdots &\vdots &\cdots &\vdots \\
\infty &*& \ldots &*
\end{array}
\right ),
$$
where some rows can consist only of infinities, and the minimal
possible value in each row is zero.

If we permute some rows of our matrix, we don't change the system of
generators of  $\tilde{\ck}(V)$ and of $\ck_x$ for each $x \in V$.
Therefore, we can assume that our matrix has the following property:
the matrix of its $Q$-orders looks like
$$
\left (
\begin{array}{ccccccccccc}
0& *&&&& \ldots &&&&&*\\
\infty &\star& \ldots &\star&0&*&&\ldots &&&*\\
\infty &\star&&& \ldots &&\star&0&*&\ldots &*\\
\vdots &&&&&\cdots &&&&\ddots &* \\
\infty &\infty &&&& \ldots &&&&&\infty \\
\vdots &&&&&\cdots &&&&& \\
\infty &\infty &&&& \ldots &&&&&\infty
\end{array}
\right ),
$$
where $\star >0$. (The rows in the bottom of matrix contain only
$\infty$.)

Clearly, the elementary transformations of rows like above lead to a
new system of generators of $\tilde{\ck}(V)$ (and of $\tilde{\ck}_x$
for each $x \in V$).

So, repeating such elementary transformations and interchanges for
the rows with zero $Q$-order,
 we will come
to a system of generators $\alpha_1, \ldots, \alpha_k$ that satisfy
the following additional property: for each $1 \le i \le k$ either
$\ord_Q(\alpha_i)=\infty$, or $\ord_Q(\alpha_i)=0$ and  $\alpha_i$
has an $l_i$-component of zero $Q$-order such that the corresponding
$l_i$-components of all other elements $\alpha_j$, $j \ne i$ has
infinite $Q$-order, see the matrix of $Q$-orders below:
$$
\left (
\begin{array}{ccccccccccccc}
0& *&\ldots &*&\infty &*& \ldots &*&\infty &*&\ldots &*\\
\infty &\star& \ldots &\star&0&*&\ldots &*&\infty &*&\ldots &*\\
\infty &\star& \ldots &\star &\infty &\star &\ldots &\star&0&*&\ldots &*\\
\vdots &\star &\cdots &\star &\infty &\star &\cdots &\star &\infty &\star &\ddots &* \\
\infty &\infty &&&&& \ldots &&&&&\infty \\
\vdots &&&&&&\cdots &&&&& \\
\infty &\infty &&&&& \ldots &&&&&\infty
\end{array}
\right ).
$$

Let $\alpha_1, \ldots , \alpha_l$ be of $Q$-order zero, and
$\alpha_{l+1}, \ldots , \alpha_k$ be of $Q$-order $\infty$. Then
\linebreak $\alpha_j=y^{k_j}\alpha''_j$, $j\ge l+1$, where
$\ord_Q(\alpha''_j)<\infty$. After multiplication of $\alpha''_j$ by
some power of $t$ we have $\alpha_j=y^{k_j} t^{m_j} \alpha'_j$,
$j\ge l+1$ for some $k_j
> 0$ and some $m_j$
such that $\ord_Q(\alpha'_j)=0$.

 We claim that the elements $\alpha_1, \ldots , \alpha_l,
\alpha'_{l+1}, \ldots ,\alpha'_k$ are generators of the $\ca
(V)$-module $\ck (V)$ such that their images generate stalks $\ck_x$
for any $x\in V$.

Indeed, if $x \in V$, $x \ne Q$, then it is clear by the choice of
elements $\alpha_1, \ldots, \alpha_k$ in the beginning, because
$\ck_x = \tilde{\ck}_x$. Now let $b\in \ck_Q$. Then $b=\sum
b_j\alpha_j$ for some $b_j\in \tilde{\ca}_Q$. We have $b_j\alpha_j
\in \ck_Q$ for all $j\ge l+1$, since $\ord_Q(b_j\alpha_j)=\infty$.
The first component of $\alpha_1$ is of zero $Q$-order, and the
$Q$-orders of the first components of all other $\alpha_i$, $i\ge 2$
are infinite. Since $b \in \ck_Q$, the $Q$-order of the first
component of $b_1\alpha_1$ must be therefore greater or equal to
zero. Hence, $\ord_Q(b_1) \ge 0$ and $b_1\in \ca_Q$. Analogously,
$b_j\in \ca_Q$ for $j\le l$. Now for $j\ge l+1$ we have
$b_j\alpha_j=b_jy^{k_j} t^{m_j} \alpha'_j$ with $k_j>0$, and
$b'_j:=b_jy^{k_j} t^{m_j} \in \ca_Q$ because $k_j >0$. So,
$$
b=\sum_{i=1}^lb_i\alpha_i+\sum_{j=l+1}^kb'_j\alpha'_j,
$$
where $b_i,b'_j\in \ca_Q$, and we are done.

Nevertheless, the sheaf $\ca$ is not weakly Noetherian. For example,
consider the following infinite increasing system of ideals in $\ca
(U)$ (for any $U \ni Q$):
$$
J_k:=\{ c=\sum_{i=l}^{\infty}c_it^i, \mbox{ where } c_i \in \cj_Q(U)
 \mbox{ and } c_i\in \cj^2_Q(U) \mbox{ for $i<-k$ }\} \mbox{.}
$$
Clearly, $J_1\subset J_2\subset \ldots$ does not stabilize.
\end{ex}

\begin{rem}
The situation described in example~\ref{ex2a} is similar to the
situation of rank $2$ valuation ring $\co^{\prime} = k[[t]] +
uk((t))[[u]]$ in $2$-dimensional local field $k((t))((u))$. The ring
$\co^{\prime}$ is also non Noetherian (see~\cite{Par}), but one can
prove that the ring $\co^{\prime}$ is coherent by the same methods
as above.
\end{rem}

\begin{ex}
\label{ex2b}
Now let's consider one more example. Let $C$ be a reduced algebraic curve over a field $k$.
Consider a ringed space $(C, \ca )$, where
$$
\ca = \{\sum_{j=N}^{\infty} \co_C\cdot t_j, \mbox{\quad} t_0=1, \mbox{\quad} t_it_j=0 \mbox{ for all $i,j\ne 0$} \}
$$
Clearly, $\ca$ is a sheaf that satisfies all conditions of definition \ref{ribbon}. So, $(C, \ca )$ is a ribbon.

Obviously, the sheaf $\ca$ is also not coherent and not weakly
Noetherian. Moreover, $\ca_0$ is not coherent. To see this, it is
enough to consider the kernel of multiplication by $t_1$. Clearly,
this kernel can not be locally of finite type.
\end{ex}

\vspace{0.5cm}

Under certain conditions on the sheaf $\ca$ of a ribbon we can prove
in the following lemma that it will be coherent, as well as any
torsion free sheaf of finite rank on this ribbon will be coherent.
(We will define torsion free sheaves later, see
definition~\ref{tfsh} and remark~\ref{remtf}).

\begin{defin}
\label{uslovie}
We will say that the sheaf $\ca$ of a ribbon $(C, \ca )$ satisfies (*), if the following condition holds:
\begin{align*}
 \mbox{ there is an affine open cover $\{U_{\alpha}\}_{\alpha\in I}$ of $C$ such that for any $\alpha \in I$ }\\
\mbox{ there is $k>0$ and an invertible section $a\in \ca_k
(U_{\alpha})\subset \ca (U_{\alpha})$.} && \hspace{2.7cm}\text{(*)}
\end{align*}
\end{defin}

\begin{defin}
For an open set $U$ we define the function of order $\ord_U$ on $\ca
(U)$ in the following way: if an element $b\in \ca_l(U)\backslash
\ca_{l+1}(U)$, then $\ord_U(b)=l$. Sometimes, if it is clear from
the context, we will omit the index $U$.
\end{defin}

Now we prove the following lemma.

\begin{lemma}
\label{lemma2} Let the sheaf $\ca$ of a ribbon $(C,\ca )$ satisfy
(*). Then it is weakly Noetherian and coherent. Moreover, for any
affine open subset $U$ of $C$ the ring $\ca (U)$ is a Noetherian
ring.
\end{lemma}

\Proof.  Let $\{U_{\alpha}\}$ be the cover from (*). For an open
$U_{\alpha}\subset C$ let $a\in \ca^*(U_{\alpha})$, $a\in
\ca_k(U_{\alpha})$, $k>0$. From the definition of ribbon
(definition~\ref{ribbon}) it follows that $a^{-1} \in
\ca_l(U_{\alpha}) \backslash \ca_{l+1}(U_{\alpha})$, where $l\le
-k$. Clearly, $\ca (U_{\alpha})=\ca_0(U_{\alpha})_a$. By
propositions~\ref{prop1} and~\ref{prop3}, the ring
$\ca_0(U_{\alpha})/\ca_{-l}(U_{\alpha})$ is Noetherian.

  Let $\tilde{I}\subset \ca (U_{\alpha})$ be an ideal. Let $I=\tilde{I}\cap \ca_0(U_{\alpha})$.
  Set $I_{-l}=I/I\cap \ca_{-l}(U_{\alpha})$. Let $\bar{g}_1, \ldots , \bar{g}_s$
  be generators of $I_{-l}$ in $\ca_0(U_{\alpha})/\ca_{-l}(U_{\alpha})$, and $g_1, \ldots ,g_s$
  be
  any their representatives in $I$.
  Let $x\in I$ be any element, $x \in \ca_j(U_{\alpha})\backslash \ca_{j+1}(U_{\alpha})$.
If $j< -l$, then there are $b_1, \ldots ,b_s\in \ca_0(U_{\alpha})$
   such that $x-\sum_m b_mg_m \in \ca_i(U_{\alpha})\backslash \ca_{i+1}(U_{\alpha})$ with $i\ge -l$.
   If $j\ge -l$, then $a^{-1}x\in I$, and for some $m \ge 1$ we have
   $a^{-m} x \in \ca_i(U_{\alpha}) \backslash \ca_{i+1}(U_{\alpha})$ with $0 \le i <-l$.
   We iterate this procedure. Since $\ord(a) > 0$, and $\ca_0(U_{\alpha})$ is a complete and
   Hausdorff space,
    we can deduce that $g_1, \ldots ,g_s$ generate $I$, hence $\tilde{I}$.
    So, $\ca (U_{\alpha})$ is a Noetherian ring.

Analogously we can show that $\ca_0(U_{\alpha})$ is also a Noetherian ring.
Namely, for an ideal $J\subset \ca_0(U_{\alpha})$ let $\tilde{J}$ be an ideal generated by $J$ in $\ca (U_{\alpha})$.
If $\tilde{J}=(1)$, then $a^r\in J$ for some $r>0$. For any $i\ge -lr$ we have $(a^r)\supseteq \ca_i(U_{\alpha})$.
Therefore, elements $g_1,\ldots ,g_s$, whose images in $\ca_0(U_{\alpha})/\ca_{-lr}(U_{\alpha})$
generate the ideal $J/J\cap \ca_{-lr}(U_{\alpha})$, and the element $a^r$ will generate the ideal $J$.

If $\tilde{J}\neq (1)$, then $\tilde{J}=(g_1, \ldots ,g_s)$ as
above, where $g_1, \ldots ,g_s\in \ca_0(U_{\alpha})$. As it was
shown above, for any sufficiently large $i$ an element $x\in J\cap
\ca_i(U_{\alpha})$ can be written as $x=a^h\sum_mb_mg_m$ with
$b_1,\ldots, b_s \in \ca_0(U_{\alpha})$, $h>0$. On the other hand,
for a sufficiently large $h$ we have $a^hg_1, \ldots, a^hg_s\in J$.
Therefore, there exists a natural $N$ such that for any $x\in J\cap
\ca_i(U_{\alpha})$ with $i>N$ we have $x\in (a^hg_1,
\ldots,a^hg_s)\subset \ca_0(U_{\alpha})$. Now, if $g'_1, \ldots,
g'_t\in J$ are representatives of generators of the ideal $J/J\cap
\ca_N(U_{\alpha})$, then the system $g'_1, \ldots ,g'_t,
a^hg_1,\ldots,a^hg_s$ is a system of generators of the ideal $J$.

To show that $\ca$ is coherent, it is enough to prove that the sheaf
$\ck$ from the definition of a coherent sheaf (see
remark~\ref{rem1}) is locally of finite type for each $U_{\alpha}$.

For any open $V\subset U_{\alpha}$ we have $\ck (V)=(\ck'(V))_a$ and $\ca (V)=(\ca_0(V))_a$, where
$$
\ck'=ker ((\ca_0 |_{U_{\alpha}})^{\oplus q}
\overset{(f_1a^{k},\ldots
,f_qa^{k})}{\longrightarrow} (\ca_0 |_{U_{\alpha}}))
$$
for sufficiently large $k$ (as in example \ref{ex2a}).
We also have
$$
\limproj\limits_{n \ge 0} \ca_0(V)/ a^n \ca_0(V)  = \ca_0(V)
\mbox{,}
$$
because for ideal $(a) = a \ca_0(V)$
 we always have $\ca_n\supseteq (a)^n\supseteq \ca_i(V)$ for $i\ge
-ln$ and for any $n$,  and $(a)^n\supseteq \ca_i(V)\supseteq (a)^i$
for $n\le [i/(-l)]$.

Combining all together, we obtain that the following locally ringed
spaces are isomorphic:
$$(U_{\alpha}, \ca_0 |_{U_{\alpha}})\simeq
\widehat{(\Spec \ca_0(U_{\alpha}))}_Y \mbox{,} $$
 where $Y$ is a
closed subscheme of $\Spec \ca_0(U_{\alpha})$ given by the ideal
$(a)$, and the formal Noetherian scheme $\widehat{(\Spec
\ca_0(U_{\alpha}))}_Y$ is a completion of the scheme $\Spec
\ca_0(U_{\alpha})$ along $Y$. So, by \cite[ch.I, \S 10.10]{EGAI},
the sheaf $\ca_0 |_{U_{\alpha}}$ is coherent, and the sheaf $\ck'$
of $\ca_0 |_{U_{\alpha}}$-modules is locally of finite type.
Therefore the sheaf $\ck$ of $\ca |_{U_{\alpha}}$-modules is locally
of finite type.

We show the last property of the lemma. At first, we note that for
any open $V\subset U_{\alpha}$ the ring $\ca (V)$ satisfies (*) and
therefore is Noetherian, as it was shown above. Since for an open
affine $U=\underline{\Spec} B$ there is a base of topology
consisting of open sets $D(f)\simeq \underline{\Spec} B_f$ and any
affine set is quasicompact, we can cover the set $U$ by finite
number of affine open sets $U_i\simeq \underline{\Spec} B_{f_i}$
such that the rings $\ca (U_i)$ satisfy (*) and are Noetherian. By
definition of a ribbon and by proposition \ref{prop3}, we can take
$B=\ca_0(U)/\ca_1(U)$, and $f_i\in \ca_0(U)/\ca_1(U)$, $f_i$
generate the ideal $(1)$ of the ring $B$.

Now we prove the following statement. Let $I\subset \ca (U)$ be an
ideal and $\phi_i:\ca (U)\rightarrow \ca (U_i)$ are the restriction
homomorphisms, $i=1, \ldots, r$. Then
$$
I=\bigcap_i \phi_i^{-1}(\phi_i(I)\cdot \ca (U_i)) \mbox{.}
$$

Obviously, we have $I\subset \bigcap_i \phi_i^{-1}(\phi_i(I)\cdot
\ca (U_i))$. Now let
 $$b\in \bigcap_i \phi_i^{-1}(\phi_i(I)\cdot \ca
(U_i)) \quad , \quad b\in \ca_k(U) \backslash \ca_{k+1}(U)
\mbox{.}$$ Let
$$\phi_i(b)=\sum_{j=1}^{r_i} \phi_i(a_j)g_j \mbox{,}$$
where $g_j\in \ca
(U_i)$, $a_j\in I$.
 We have $\phi_i(b)\in \ca_k(U_i)$ and therefore
$$
\phi_i(b) \mod  \ca_{k+1}(U_i) \; = \; \sum_{j=1}^{r_i} (\phi_i(a_j)
\mod \ca_{k+1-\ord (g_j)}(U_i))(g_j \mod  \ca_{k+1-\ord (a_j)}(U_i))
\mbox{.}
$$
We consider the homomorphisms
$$\bar{\phi}_i^j \; : \; \ca_{\ord (g_j)}(U)/\ca_{k+1-\ord (a_j)}(U)
\longrightarrow  \ca_{\ord (g_j)}(U_i)/\ca_{k+1-\ord (a_j)}(U_i)
\mbox{,}
$$
which are induced by $\phi_i$. By proposition \ref{prop1}, the sheaf
$\ca_{\ord (g_j)}/\ca_{k+1-\ord (a_j)}$ is a coherent sheaf on the
scheme $X_l=(C, \ca_0/\ca_{l+1})$, where $l=k-\ord (a_j)-\ord (g_j)$
(we assume that $l\ge0$, since otherwise our sheaf is trivial and
there is nothing to prove). Therefore, $\bar{\phi}_i^j$ is a
localization map, and for any element
 $x\in \ca_{\ord (g_j)}(U_i)/\ca_{k+1-\ord (a_j)}(U_i)$ there exists
  a natural $n$ such that $f_{ij}^nx=\bar{\phi}_i^j(\tilde{x})$,
  where $\tilde{x}\in \ca_{\ord (g_j)}(U)/\ca_{k+1-\ord (a_j)}(U)$,
  \linebreak
   $f_{ij}\in \ca_0(U)/\ca_{l+1}(U)$ and $f_{ij} \mod \ca_1(U)=f_i$
   (see \cite[lemma 5.3]{Ha}).
    Note that we can choose $f_{ij}=\tilde{f}_i \mod \ca_{l+1}(U)$,
    where $\tilde{f}_i$ is a fixed representative of $f_i$ in $\ca_0(U)$, for all $j$.
    Hence there exists a natural $N$ such that
$$
\phi_i(\tilde{f}_i^N)\phi_i(b) \mod  \ca_{k+1}(U_i) \; = \;
\sum_{j=1}^{r_i} (\phi_i(a_jg'_j) \mod  \ca_{k+1}(U_i)) \mbox{,}
$$
where $g'_j\in \ca (U)$ and
$$\phi_i(g'_j) \mod  \ca_{k+1-\ord
(a_j)}(U_i) \; = \; \phi_i(\tilde{f}_i^N)g_j \mod  \ca_{k+1-\ord
(a_j)}(U_i) \mbox{.} $$ Let $k'$ be an integer such that $a_jg'_j\in
\ca_{k'}(U)$ for any $j$ (note that $k'\le k$). Then, repeating the
arguments above to the coherent sheaf $\ca_{k'}/\ca_{k+1}$ we obtain
that there exists a natural $M$ such that
$$
\tilde{f}_i^M b \mod  \ca_{k+1}(U)\; \in \; I \mod  \ca_{k+1}(U).
$$
Note that we can choose $M$ unique for all $i$ and that the elements
$\tilde{f}_i^M$ generate the ideal $(1)$ in $\ca_0(U)$, i.e. $\sum_i
c_i\tilde{f}_i^M=1$ for some $c_i \in \ca_0(U)$. Therefore,
$$
b \mod  \ca_{k+1}(U) \; = \; \sum c_i\tilde{f}_i^Mb \mod
\ca_{k+1}(U) \; \in \; I \mod  \ca_{k+1}(U).
$$
So, there exists $b_1\in I$, $b_1\in \ca_k(U)$ such that
 $(b-b_1) \in \bigcap_i \phi_i^{-1}(\phi_i(I)\cdot \ca (U_i))$ and $\ord (b-b_1)>\ord (b)$.
 We repeat the arguments above for the element $b-b_1$ and so on.
 Since the ring $\ca (U)$ has complete and Hausdorff topology,
 we obtain that $b\in I$.

Now it is easy to show that the ring $\ca (U)$ is Noetherian. Let
$I_1\subset I_2\subset \ldots $ be an increasing chain of ideals in
$\ca (U)$. Then for each $i$ the chain
$$
\phi_i(I_1)\cdot \ca (U_i)\subset \phi_i(I_2)\cdot \ca (U_i)\subset \ldots
$$
is stable, since $\ca (U_i)$ is a Noetherian ring. Since there are
only finite number of $i$, the first chain is also stable, where
from $\ca (U)$ is a Noetherian ring.
\begin{flushright}
$\square$
\end{flushright}

\begin{defin}
A ribbon $(C, \ca)$ is called  algebraizable if it is locally
isomorphic on $C$ to a ribbon from example~\ref{ex1}.
\end{defin}

\begin{ex}
The sheaf $\ca =\co_{\hat{X}_C}(*C)$ with the filtration
$\ca_i=\co_{\hat{X}_C} (-iC)$ on a surface $X$ with an effective
reduced Cartier divisor $C$ from example \ref{ex1} satisfies the
conditions of lemma~\ref{lemma2}. Indeed, the local equation of $C$
in $X$ is an invertible element that belong to $\ca_1(U)$, and its
inverse belongs to $\ca_{-1}(U)$.

In particular, it follows that the ribbons from example~\ref{ex2},
example~\ref{ex2a} and example~\ref{ex2b} are not algebraizable,
because they are not weakly Noetheriean.
\end{ex}

\begin{rem}
Structure sheaves of algebraizable ribbons satisfy  more pretty
property, which is useful in studying of the Picard group of a
ribbon, see proposition~\ref{lemma5} below.
\end{rem}

\begin{ex}
\label{ex2d} We consider an example of a ribbon with weakly
Noetherian and coherent structure sheaf $\ca$, but which is not
algebraizable. It can be constructed in the same way as in example
\ref{ex2b}.

Let $C$ be a reduced algebraic curve over a field $k$. Consider a ringed space $(C, \ca )$, where
$$
\ca = \{\sum_{j=N}^{\infty} \co_C\cdot t_j, \mbox{\quad} t_0=1, \mbox{\quad}
t_{2i}=t_2^i, \mbox{\quad} t_{2i+1}=t_1t_2^i, \mbox{\quad}
 t_1^2=0  \}.
$$
Clearly, $\ca$ is a sheaf that satisfies all conditions of
definition \ref{ribbon}. So, $(C, \ca )$ is a ribbon. By lemma
\ref{lemma2} $\ca$ is a weakly Noetherian and coherent sheaf (since
$t_2$ is an invertible section of $\ca (U)$ for any open $U\subset
C$). But $(C,\ca )$ is not algebraizable, because if it were
algebraizable, there should be an open affine cover of $C$ such that
for any open $U$ from this cover there exists an invertible element
$a$ that belong to $\ca_1(U)\backslash \ca_2(U)$ and $a^{-1} \in
\ca_{-1}(U)\backslash \ca_0(U)$. Obviously, there are no such
sections in $\ca (U)$ for any $U$.
\end{ex}

\subsection{Analytic ribbons}
 When a ground field is $\bf C$, we can also work in
the analytic category to define ribbons over $\bf C$, replacing
"algebraic coherent" by "analytic coherent sheaf" (for
$\ca_i/\ca_{i+1}$, $i \in \dz$) in definition~\ref{ribbon}. Then we
obtain the notion of an analytic ribbon $(C, \ca)$.

We define an analytic ind-pro-coherent sheaf $\cf$ on analytic
ribbon  $\xo_{\infty} = (C, \ca)$ as a filtered sheaf of
$\ca$-modules (with a descending filtration by subsheaves)
satisfying properties \ref{s1}, \ref{s3}, \ref{s4} and the property
$$
\mbox{  {\ref{s2}$'$}.  $\cf_j/\cf_{j+1}$ is an analytic coherent
sheaf on $C$  for any $j \in \dz$}$$
  instead of property~\ref{s2} of definition~\ref{ipqs}.

\begin{rem}
Since the underlying topological space is non-Noetherian in this
case, we have to take the sheaf $\cf$ associated with the presheaf
$\cf': V \mapsto \liminj\limits_i \cf_i(V)$.
\end{rem}

We have the following proposition (compare with
proposition~\ref{prop3}).
\begin{prop}
\label{prop3'} We have the following properties for an analytic
ind-pro-coherent sheaf $\cf$ on an analytic ribbon $\xo_{\infty} = (C,
\ca)$,  where $C$ is an irreducible complex algebraic curve.
\begin{enumerate}
\item
\label{sta1} $\cf_i/\cf_{i+j+1}$ is an analytic coherent sheaf on
$X_j$ for any $j \ge 0$, $i \in \dz$.
\item
\label{sta2}
$\cf_i(U)/\cf_j(U) \rightarrow (\cf_i/\cf_j)(U)$ is an
isomorphism for $i<j$ and for Stein open sets $U\subset C$.
\item
\label{sta3} $H^q(U,\cf_i)=0$ for any Stein open subset $U \subset
C$ and $q>0$, $i \in \dz$.
\end{enumerate}
\end{prop}
\begin{rem} \label{rema6}
We note that every complex analytic space of dimension $1$, which has no
compact irreducible components, is a Stein space (see, for
example,~\cite{N}).
\end{rem}
\Proof. The proof of statement~\ref{sta1} and statement~\ref{sta2}
of this proposition is the same as in proposition \ref{prop3}. (We
use that for any analytic coherent sheaf $\cg$ on a Stein space $U$
we have $H^q(U, \cg)= 0$ for $q > 0$.)

Now we prove statement~\ref{sta3} of the proposition. By
remark~\ref{rema6} we have that any open subset $V$ of a Stein
subset $U \subset C$ is a Stein space again. Therefore, if
$\{U_{\alpha}\}$ is an open covering of $U$, then every open
$U_{\alpha}$ is a Stein space. Let
$\check{\cc}^{\bullet}(\{U_{\alpha}\}, \cf_i)$ be the
\v{C}ech-complex of this covering for the sheaf $\cf_i$. We obtain
that
$$
\check{\cc}^{\bullet}(\{U_{\alpha}\},
\cf_i)=\limproj\limits_{j>i}\check{\cc}^{\bullet}(\{U_{\alpha}\},
\cf_i/\cf_j)\mbox{.}
$$
We consider the following natural complex $D^{\bullet}_i$ for any $i
\in \dz$:
$$
0 \lto \cf_i(U) \lto \check{\cc}^{0}(\{U_{\alpha}\}, \cf_i) \lto
\check{\cc}^{1}(\{U_{\alpha}\}, \cf_i) \lto \ldots \mbox{,}
$$
i.e. $D^{n}_i = 0$ for $n < -1$, $D^{-1}_i = \cf_i(U)  $, and
$D^n_i= \check{\cc}^{n}(\{U_{\alpha}\}, \cf_i)$ for $n \ge 0$.

We have that for any $i \in \dz$
$$
D_i^{\bullet} = \limproj\limits_{j>i} D_{i,j}^{\bullet} \mbox{,}
$$
where  the complex $D_{i,j}^{\bullet}$ is defined in the following
natural way for any $j \ge i \in \dz$: $D^{n}_{i,j} = 0$ for $n <
-1$, $D^{-1}_{i,j} = (\cf_i / \cf_j)(U)  $, and
$D^n_{i,j}=\check{\cc}^{n}(\{U_{\alpha}\}, \cf_i / \cf_j)$ for $n
\ge 0$.

From statement~\ref{sta2} of this proposition we have that for any
fixed $i \in \dz$,  for any $n \in \dz$ the projective system
$(D^n_{i,j}, j \ge i)$ satisfies the ML-condition, because the maps
in this projective system are surjective maps.

For any $j \ge i \in \dz$ the complex $D^{\bullet}_{i,j}$ is an
acyclic complex, because the \v{C}ech cohomology
$$\check{H}^0
(\{U_{\alpha} \}, \cf_i/ \cf_j) = (\cf_i /\cf_j) (U)  \mbox{,}
$$
$$
\check{H}^n (\{ U_{\alpha} \}, \cf_i / \cf_j) = H^n (U, \cf_i /
\cf_j) = 0 \quad  \mbox{for any} \quad n
>0 \mbox{.}
$$
Therefore for any $i \in
\dz$ the complex $D^{\bullet}_i$ is an acyclic complex, as it
follows from the following lemma.
\begin{lemma} \label{lemma*}
Let $(K_l^{\bullet}, l \ge 0 )$ be a projective system of acyclic
complexes $K_l^{\bullet}$ of abelian groups.  We suppose  that for
any $n \in \dz$ the projective system $(K_l^n, l \ge 0)$ satisfies
the ML-condition. Then the complex
$$
K^{\bullet} = \limproj\limits_{l \ge 0} K_l^{\bullet}
$$
is an acyclic complex.
\end{lemma}
\Proof \: of the lemma. Let  maps $d^n_{l} : K_l^n  \to K_l^{n+1}$,
$n \in \dz$ be the differentials in complex $K_l^{\bullet}$, $l \ge
0$. We have the following exact sequences:
\begin{equation} \label{compl}
0 \lto \Ker d^n_l \lto K_l^n \lto \I d^{n}_l \lto 0 \mbox{.}
\end{equation}
Since the complex $K^{\bullet}_l$ is an acyclic complex, we have
that $\I d^{n}_l = \Ker d^{n+1}_l$ for any $n$.

Since for any $n$ the projective system $(K_l^n, l \ge 0 )$
satisfies ML-condition, we obtain from exact sequence~(\ref{compl})
that for any $n$ the projective system $(\I d^{n-1}_l, l\ge 0)
=(\Ker d^{n+1}, l \ge 0)$ satisfies ML-condition. Let maps $d^n :
K^n \to K^{n+1}$ be the differentials in complex $K^{\bullet}$. Now,
using lemma~\ref{lemma1} and that always $\Ker d_n =
\limproj\limits_{n \ge 0} \Ker d_n^l$ for any $n \in \dz$, we obtain
that the projective limit with respect to $l \ge 0$ of
sequences~(\ref{compl}) will give the following exact sequence for
any $n \in \dz$:
$$
0 \lto \Ker d^n \lto K^n \lto \I d^{n} \lto 0 \mbox{.}
$$

Therefore, for any $n \in \dz$ we have
$$
\I d^{n} = \limproj\limits_{l \ge 0} \I d^{n}_l = \limproj\limits_{l
\ge 0} \Ker d^{n+1}_l = \Ker d^{n+1} \mbox{.}
$$
Therefore the complex $K^{\bullet}$ is an acyclic complex. The lemma
is proved.
\begin{flushright}
$\square$
\end{flushright}
Now we finish the proof of  proposition~\ref{prop3'}. We have proved that
$\check{H}^q(\{U_{\alpha}\},\cf_i)=0$ for any $i \in \dz$ and any $q
> 0$. Therefore $\check{H}^q(U,\cf_i)=0$ for any $i \in \dz$ and any $q > 0$.
Hence, for any $i \in \dz$
$$
H^1(U,\cf_i)= \check{H}^1(U,\cf_i)=0 \mbox{.}
$$

 Furthermore, we have  a spectral sequence with initial
term
\begin{equation} \label{specseq}
E_2^{pq}=\check{H}^p(U,\ch^q(\cf_i ))\Longrightarrow H^{p+q}(U,\cf_i
) \mbox{,}
\end{equation}
where $\ch^q(\cf_i )$ is the presheaf $V\subset U \mapsto
\ch^q(\cf_i )(V)=H^q(V,\cf_i )$ (see \cite{Gr}). So, since any open
subset $V \subset U$ is  a Stein subspace again, we have in our
situation
$$
H^2(U,\cf_i )=\check{H}^0(U,\ch^2(\cf_i )) \mbox{.}
$$
To obtain that $H^2(U,\cf_i )=0$ it is sufficient to show that for
any  point $x\in C$
$$\liminj\limits_{x\in V\subset C }H^2(V,\cf_i)=0
\mbox{.}$$ It follows from the following fact (\cite[lemma
3.8.2]{Gr}: for any point $x \in C$, for any $p >0$
\begin{equation} \label{equality}
\liminj\limits_{x\in V\subset C }H^p(V,\cf_i)=0
\end{equation}

 Now by induction on $q$, by the same methods as for $q=2$, using spectral sequence~(\ref{specseq})
 and equality~(\ref{equality}), we
obtain that $H^q(U,\cf_i )=0$ for all $q>0$. The proposition is
proved.
\begin{flushright}
$\square$
\end{flushright}

\begin{corol}
Let $\xo_{\infty} = (C, \ca)$ be an analytic ribbon. Let $\cf$ be an
analytic ind-pro-coherent sheaf on $\xo_{\infty}$, and $C$ be an irreducible
compact space.
\begin{enumerate}
\item
If $C = U_1 \cup U_2$, where  $U_1$ and $U_2$ are Stein open
subsets, then we have an exact sequence for any $i \in \dz$
$$
0 \rightarrow H^0(C,\cf_i ) \rightarrow H^0(U_1,\cf_i )\oplus
H^0(U_2,\cf_i ) \rightarrow H^0(U_1 \cap U_2, \cf_i ) \rightarrow
H^1(C,\cf_i )\rightarrow 0 \mbox{.}
$$
\item
$H^*(C,\cf_i )=\limproj\limits_{j > i}H^*(C, \cf_i/\cf_{j}) $, $i
\in \dz$.
\item
$ H^q(C, \cf_i) =0$ for $q > 1$, $i \in \dz$.
\end{enumerate}
\end{corol}
\Proof \: is similar to the proof of corollary~\ref{col1} of
proposition~\ref{prop3}. We have to use the following Mayer-Vietoris
exact sequence for a sheaf $\cg$ on $C$:
$$
\ldots \lto H^{k-1}(U_1 \cap U_2, \cg) \lto H^k(C, \cg) \lto H^k
(U_1, \cg) \oplus H^k(U_2, \cg) \lto \ldots
$$
\begin{flushright}
$\square$
\end{flushright}

\section{The Picard group of a ribbon}
\label{Pic}

We recall that for a ringed space $\xo_{\infty}=(C,\ca )$ the Picard
group $ Pic (\xo_{\infty})= H^1(C, \ca^*)$, and  for the ringed
space $X_{\infty}= (C, \ca_0) $  also the Picard group $Pic
(X_{\infty}) = H^1 (C, \ca_0^*)$.

\begin{prop} \label{propos5}
Let $\xo_{\infty}=(C,\ca )$ be a ribbon over an Artinian ring $A$.
We suppose that $C$ is either  projective, or affine curve over
$\Spec A$. Then
$$
Pic(X_{\infty})=\limproj\limits_{i \ge 0} Pic(X_i) \mbox{.}
$$
\end{prop}
\Proof. We denote for any $j \ge i \ge 0$ the following sheaves
$\cg_{i,j} = \frac{1+\ca_{i+1}}{1+\ca_{j+1}}$ on $C$. Then we have
the following exact sequences:
\begin{equation} \label{equ0}
1 \lto \cg_{i,j} \lto  \co_{X_j}^*  \lto \co_{X_i}^*  \lto 1
\mbox{.}
\end{equation}

We denote for any $i \ge 0$ the following sheaf $\cg_i = 1
+\ca_{i+1} \subset \ca_0^*$ on $C$. Then we have
$$
\cg_i = \limproj\limits_{j \ge i} \cg_{i,j} \mbox{.}
$$

For any $j \ge i \ge 0$ we have the following exact sequences:
\begin{equation} \label{sss1}
1 \lto \cg_{j, j+1} \lto \cg_{i, j+1} \lto \cg_{i, j} \lto 1
\mbox{.}
\end{equation}
For any $j \ge 0$ we have $\cg_{j, j+1} \simeq \ca_{j+1}/
\ca_{j+2}$. Therefore from sequence~(\ref{sss1}) we obtain that for
any affine open subset $U \subset C$ the maps $H^0(U, \cg_{i,j+1})
\to H^0(U, \cg_{i,j})$ are surjective for any $j \ge i \ge 0$, by
induction on $j$ we obtain that $H^1(U, \cg_{i,j})=0$ for any $j \ge
i \ge 0$. Therefore, arguing as in the proof of assertion~\ref{i3}
of proposition~\ref{prop3}, we obtain that $H^1(U, \cg_i) =0$ for
any $i \ge 1$.

Since $C$ is a curve over an Artinian ring, there are some affine
open subsets $U_1$ and $U_2$ of $C$ such that $C= U_1 \cup U_2$.
Therefore the following Mayer-Vietoris sequence is exact:
\begin{equation} \label{equ1}
0\rightarrow H^0(C,\cg_i ) \rightarrow H^0(U_1,\cg_i )\oplus
H^0(U_2, \cg_i )\rightarrow H^0(U_1 \cap U_2, \cg_i ) \rightarrow
H^1(C,\cg_i )\rightarrow 0 \mbox{.}
\end{equation}

Also for any $j \ge  i \ge 0$ we have the following exact sequences:
\begin{equation} \label{equ2}
0\rightarrow H^0(C,\cg_{i,j} ) \rightarrow H^0(U_1,\cg_{i,j} )\oplus
H^0(U_2, \cg_{i,j} )\rightarrow H^0(U_1 \cap U_2, \cg_{i,j})
\rightarrow H^1(C,\cg_{i,j} )\rightarrow 0 \mbox{.}
\end{equation}

We note that for any fixed $i \ge 0$ the projective system
$(H^0(U_1,\cg_{i,j} ) \oplus H^0(U_2, \cg_{i,j}), j \ge i)$
satisfies the ML-condition, because the maps in this system are
surjective. By the same reason, if the curve $C$ is affine, then for
any fixed $i \ge 0$ the projective system $(H^0(C,\cg_{i,j}), j \ge
i)$ satisfies the ML-condition. If the curve $C$ is projective, then
we consider the following exact sequences which follow from
sequences~(\ref{equ0}):
\begin{equation} \label{equ4}
0 \lto H^0(C, \cg_{i,j}) \lto H^0(C, \co^*_{X_j}) \lto H^0(C,
\co^*_{X_i})   \mbox{.}
\end{equation}
We have that $A$-modules $(H^0(C, \co_{X_j}), j \ge 0)$ satisfy the
ML-condition, because they are Artinian $A$-modules. Therefore the
groups $H^0(C, \co^*_{X_j}) = H^0(C, \co_{X_j})^*$ satisfy the
ML-condition as invertible elements of the corresponding algebras
for which: 1) we have ML-condition and 2) maps in projective system
have nilpotent kernels. Whence, for the fixed $i \ge 0$ from exact
sequence~(\ref{equ4}) we obtain that the projective system $(H^0(C,
\cg_{i,j}), j \ge i)$ satisfies the ML-condition as the kernels of
the maps to  the constant group $H^0(C, \co^*_{X_i})$.

Now we apply lemma~\ref{lemma1} to obtain that for the fixed $i \ge
0$ exact sequence~(\ref{equ1}) is the projective limit of exact
sequences~(\ref{equ2}) with respect to $j \ge i$. Therefore we have
$$
H^1(C, \cg_i) = \limproj \limits_{j \ge i} H^1 (C, \cg_{i,j})
\mbox{.}
$$

Let $i =0$, then from exact sequence~(\ref{equ0}) we obtain the
following exact sequence for $j \ge 0$:
$$
0 \lto H^0 (C, \cg_{0,j}) \lto H^0(C, \co_{X_j}^*) \lto H^0(C,
\co_C^*) \lto \qquad \qquad \qquad \qquad \qquad \qquad
$$
$$
\qquad \qquad \qquad \qquad \qquad \qquad
 \lto H^1(C, \cg_{0,j})
\lto H^1(C, \co_{X_j}^*) \lto H^1(C, \co_C^*) \lto 0 \mbox{.}
$$
Every term of this sequence satisfies the ML-condition with respect
to $j$. (For the zero cohomology it is proved above in this proof,
for the first cohomology it follows from the absence of the second
cohomology on the curve. ) Therefore, using lemma~\ref{lemma*}, we
obtain that the following sequence is exact:
$$
0 \to \limproj \limits_{j \ge 0} H^0 (C, \cg_{0,j}) \to
\limproj\limits_{j \ge 0} H^0(C, \co_{X_j}^*) \to H^0(C, \co_C^*)
\to \qquad \qquad \qquad \qquad \qquad \qquad
$$
\begin{equation} \label{eeee1}
\qquad \qquad \qquad \qquad \qquad \qquad
 \to \limproj \limits_{j \ge 0}H^1(C, \cg_{0,j})
\to \limproj \limits_{j \ge 0}H^1(C, \co_{X_j}^*) \to H^1(C,
\co_C^*) \to 0 \mbox{.}
\end{equation}

From exact sequence
$$
1 \lto \cg_0 \lto \ca_0^*  \lto \co_C^* \lto 1
$$
we obtain the following exact sequence:
$$
0 \lto H^0 (C, \cg_{0}) \lto H^0(C, \ca_0^*) \lto H^0(C, \co_C^*)
\lto \qquad \qquad \qquad \qquad \qquad \qquad
$$
\begin{equation} \label{eeee2}
\qquad \qquad \qquad \qquad \qquad \qquad
 \lto H^1(C, \cg_{0})
\lto H^1(C, \ca_0^*) \lto H^1(C, \co_C^*) \lto 0 \mbox{.}
\end{equation}

We have the natural map of exact sequence~(\ref{eeee2}) to exact
sequence~(\ref{eeee1}), and we know that the maps at each term
except for  one term of the sequence are isomorphisms. But then it
follows that the residuary map
$$
H^1(C, \ca_0^*)  \lto \limproj \limits_{j \ge 0}H^1(C, \co_{X_j}^*)
$$
is also an isomorphism.
\begin{flushright}
$\square$
\end{flushright}

\begin{corol}
Under conditions of proposition~\ref{propos5} we suppose that $C$ is
an affine curve. Then $Pic(X_{\infty}) = Pic(C)$.
\end{corol}
\Proof\: follows from the proposition and from $H^1(C, \cg_{j, j+1})
=1$ for any $j \ge 0$.
\begin{flushright}
$\square$
\end{flushright}

Let $\xo_{\infty} = (C, \ca)$ be a ribbon over the field $k$. For a
point $x \in C$ we denote by $\ca_{0,x}$ the local ring which is a
stalk of the sheaf $\ca_0$ at the point $x$. Let $\cm_x$ be the
maximal ideal of $\ca_{0,x}$. Further we will need to compare the
following  two rings.
\begin{defin} \label{defrin}
 We denote by $\widehat{\ca}_{0,x}$ the $\cm_x$-adic completion of the ring
 $\ca_{0,x}$.
 We denote by $\tilde{\ca}_{0,x}$ the ring
 $$\tilde{\ca}_{0,x} = \limproj\limits_i\widehat{\co}_{X_i,x} \mbox{.}$$
\end{defin}

\begin{prop}
\label{prop5}
\begin{enumerate}
\item \label{as1pr7}
We have the following commutative diagram of morphisms of local
rings
$$
\begin{array}{ccc}
\ca_{0,x}&\rightarrow &\widehat{\ca}_{0,x}\\
\parallel && \downarrow \lefteqn{\alpha} \\
\ca_{0,x}&\rightarrow & \tilde{\ca}_{0,x}
\end{array}
$$
where the horizontal arrows are injective.
\item \label{as2pr7} If $\dim_k
(\cm_x/\cm_x^2)<\infty $, then the ring $\widehat{\ca}_{0,x}$ is
Noetherian and $\alpha$ is surjective, and the Krull dimension of $
\tilde{\ca}_{0,x}$: $dim \tilde{\ca}_{0,x} \ge 2$. Furthermore,
$\tilde{I}_j=I_j \widetilde{\ca}_{0,x}$, where $\tilde{I}_j=Ker
(\widetilde{\ca}_{0,x} \rightarrow \widehat{\co}_{X_i,x})$,
$I_j=\ca_{j,x}$.
\end{enumerate}
\end{prop}

\Proof. We prove assertion~\ref{as1pr7} of the proposition. We
define a linear topology on \linebreak $A:=\ca_{0,x}$ by taking as
open ideals all ideals $\cq$ of finite colength which contain some
ideal $I_i:=\ca_{i,x}$. Thus the set $\{\cq\}$ of ideals contains
the ideals $I_i+\cm_x^n$ for all $i,n$, since $A/I_i$ is Noetherian,
and so it is coarser or equivalent to the $\cm_x$-adic topology, and
it is separated (since for any $a\ne 0$ in $A$ there is $I_i$ with
$a\neq 0 \mbox{ mod } I_i$, and $n$ with $a \mbox{ mod } I_i\notin
\cm_x^n(A/I_i)$, hence $a\notin \cm_x^n+I_i=\cq$).

Hence assertion~\ref{as1pr7} follows, since $\tilde{A}$ is the
completion of $A$ with respect to the $\{\cq\}$-topology.

We prove assertion~\ref{as2pr7} of the proposition. We recall that
$k(x) = A/\cm_x$. If \linebreak $\dim_{k(x)}
(\cm_x/\cm_x^2)=n<\infty $, then $gr_{\cm_x}(A)$ is Noetherian (as
an image of a surjection $Sym_{A/\cm_x}(\cm_x/\cm_x^2)\rightarrow
gr_{\cm_x}(A)$) and $\dim_{k(x)} (A/\cm_x^{k+1})\le C_{n+k}^n$.
Therefore $\widehat{A}$ is Noetherian by \cite[ch. III, \S 2.9,
corol.2]{Bur} (since
$gr_{\{\widehat{\cm^i_x}\}}(\widehat{A})=gr_{\cm}(A)$), and
$\widehat{A}$, $\tilde{A}$ both carry a linear topology, which is
linearly compact. (The topology of $\widehat{A}$ is linearly
compact, since $\dim_k \widehat{A} / \widehat{\cm_x^i} < \infty$,
see \cite[ch.III, \S 2]{Bur}). Since $\alpha$ is a continuous
homomorphism and $A$ is dence in $\tilde{A}$ and in $\widehat{A}$,
it follows that $\alpha$ is surjective with kernel
$\cap_{\cq}\widehat{\cq}=\cap_j\widehat{I}_j$ ($\widehat{I}$ is the
$\cm_x$-adic completion of an ideal). The topology of $\tilde{A}$ is
the $\tilde{\cm_x}$-adic topology.

Now we prove $\tilde{I}_j=I_j\tilde{A}$. We have
$\tilde{I}_j/\tilde{I}_{j+k}=I_j\widehat{A/I_{j+k}}=I_j\tilde{A}/\tilde{I}_{j+k}$,
hence $\tilde{I}_j=I_j\tilde{A}+\tilde{I}_{j+k}$ for any $k>0$. But
since $\tilde{A}$ is Noetherian, the ideal $I_j\tilde{A}$ (as any
ideal in $\tilde{A}$) is closed in the $\cm_x$-adic topology and the
$\{\tilde{I}_j\}$-topology is finer (since $\tilde{A}$ is linearly
compact). Therefore $\tilde{I}_j=I_j\tilde{A}$.

To prove that $dim \tilde{A}\ge 2$, we choose $u\in \cm_x$ which
lifts a nonzero divisor of $\cm_{C,x} =\cm_x/I_1$. We'll prove
$l(A/(I_{j+1}+uA)) \ge j+1$. It follows by induction on $j$ from the
following exact sequence
$$
\begin{array}{ccccccc}
0\rightarrow &I_j/I_{j+1}&\rightarrow &A/I_{j+1}& {\rightarrow} &A/I_j&\rightarrow 0\\
&\downarrow \lefteqn{u}&&\downarrow \lefteqn{u} &&\downarrow \lefteqn{u}& \\
0\rightarrow & I_j/I_{j+1} & \rightarrow & A/I_{j+1} & \rightarrow &
A/I_j & \rightarrow 0
\end{array}
$$
since we have statement~(\ref{iii}) of definition~\ref{ribbon}, so
$l((I_j/I_{j+1})/u(I_j/I_{j+1}))\ge 1$, \linebreak
$l((A/I_1)/u(A/I_1))=l(A/(I_1+uA)) \ge 1$, $
l((A/I_{j+1})/u(A/I_{j+1}))= l((A/I_j)/u(A/I_j)) +
l((I_j/I_{j+1})/u(I_j/I_{j+1}))$.

Thus $l(\tilde{A}/(\tilde{I}_j+u\tilde{A})) \ge j$, and
$l(\tilde{A}/u\tilde{A})=\infty$. Since $u$ is not a zero-divisor in
$\tilde{A}=\limproj \widehat{A/I_j}$, it follows that $dim
\tilde{A}>1$.
\begin{flushright}
$\square$
\end{flushright}

\begin{corol} \label{coroll}
We suppose that $\dim_{k(x)} (\cm_x/\cm_x^2)=2$.  Then we have the
following properties under notations of proposition~\ref{prop5}.
\begin{enumerate}
\item $\alpha$ is an isomorphism.
\item $\tilde{\ca}_{0,x}$ is a 2-dimensional regular ring.
\end{enumerate}
\end{corol}

\Proof. We  always have $dim \hat{\ca}_{0,x}\ge dim
\tilde{\ca}_{0,x}$. By \cite[ch.III, \S 3, prop.3]{Bur},  the
filtration $\{\widehat{\cm_x^i}\}$ is $\widehat{\cm_x}$-stable in
the ring $\hat{\ca}_{0,x}$. Then by \cite[prop.11.4, th. 11.14]{AM},
 we have $dim \hat{\ca}_{0,x}=\deg
\chi_{\cm_x}(n)=\deg g(n)$, where $\chi_{\cm_x}(n)$, $g(n)$ are
characteristic polynomials for the filtrations
$\{\widehat{\cm_x}^i\}$, $\{\widehat{\cm_x^i}\}$, and $2=dim \,
Sym_{A/\cm_x}(\cm_x/\cm_x^2)\ge \deg g(n)$ (since $g(n)\le
\chi_{\nu}(n)$ for all $n\gg 0$, where $\chi_{\nu}$ is the
characteristic polynomial of the ring
$(Sym_{A/\cm_x}(\cm_x/\cm_x^2))_{\nu}$, where the prime ideal
$\nu=\bigoplus_{n=1}^{\infty} S^n(\cm_x/\cm_x^2)$).

Therefore, using assertion~\ref{as2pr7} of proposition~\ref{prop5},
we have $ dim \tilde{\ca}_{0,x}=dim \hat{\ca}_{0,x}=2$ and
$\tilde{\ca}_{0,x}$ is a 2-dimensional regular ring with a prime
ideal $(0)$. Therefore, $ker (\alpha)$ must be a prime ideal, hence
$ker (\alpha)=0$,
 since otherwise $dim \hat{\ca}_{0,x}>2$.
\begin{flushright}
$\square$
\end{flushright}

\begin{prop}
\label{prop6} The group $\ca_x^*/\ca_{0,x}^*$ is non-trivial if and
only if there exists an integer $i>0$ such that
$\ca_{i,x}\ca_{-i,x}=\ca_{0,x}$. In this case the following
properties are satisfied.
\begin{enumerate}
\item
\label{18)}  All $\ca_{j,x}$ ($j\in \dz$) are finitely generated
$\ca_{0,x}$-modules.
\item
\label{28)}  The invertible sets $\ca_{j,x}$'s, i.e. those for which
$\ca_{j,x}\ca_{-j,x} =\ca_{0,x}$ form a cyclic group $\{\ca_{id,x} |
i\in\dz\}$ with some $d>0$.
\end{enumerate}
\end{prop}

\Proof. If $\ca_{i,x}\ca_{-i,x}=\ca_{0,x}$, there are finitely many
elements $a_i\in \ca_{i,x}$, $b_i\in \ca_{-i,x}$, $i=1,\ldots ,r$
such that $\sum a_ib_i=1$.

Since $\ca_{0,x}$ is a local ring there is one pair $(a_i,b_i)$ with $a_ib_i\in \ca_{0,x}^*$,
and so there is a pair $(a,b)$, $a\in \ca_{i,x}$, $b\in\ca_{-i,x}$ with $ab=1$.

Now from $ab=1$, $a\in \ca_{i,x}$, $b\in \ca_{-i,x}$, we obtain
$\ca_{i,x}= \ca_{0,x}a$, $\ca_{-i,x}=\ca_{0,x}b$. For, if $a'\in \ca_{j,x}$ and $a'b=f\in \ca_{0,x}$,
then $0=((a'-fa)b)a=a'-fa$, hence $a'=fa$.
Similarly,
$\ca_{ki,x}=\ca_{0,x}a^k$, $\ca_{-ki,x}=\ca_{0,x}b^k$, since
$a^kb^k=1$.

If $\ca_{i,x}\ca_{-i,x}=\ca_{0,x}$, $\ca_{j,x}\ca_{-j,x}=\ca_{0,x}$ and $d=gcd(i,j)$,
then $\ca_{d,x}\ca_{-d,x}=\ca_{0,x}$. For, if $\ca_{i,x}=\ca_{0,x}a$, $\ca_{j,x}=\ca_{0,x}a'$
and $d=mi+nj$, then $a^m{a'}^n\in \ca_{d,x}$, and if $b=a^{-1}$, $b'={a'}^{-1}$ $b^m{b'}^n\in \ca_{-d,x}$
and $(a^m{a'}^n)(b^m{b'}^n)=1$.

Thus, assertion~\ref{28)} of this proposition is proved. To prove
assertion~\ref{18)} of the proposition, we observe that for any
$\ca_{j,x}$ there is a multiple $k=dj$ of $d$ such that
$\ca_{k,x}\subset \ca_{j,x}$, and $\ca_{j,x}/\ca_{k,x}$ is a
finitely generated $\ca_{0,x}$-module.
\begin{flushright}
$\square$
\end{flushright}

Now we want to discuss the group $H^1(C,\ca^*)$.

\begin{prop}
\label{prop4} Let $\xo_{\infty}$ be a ribbon with an irreducible
underlying curve $C$. We assume that the function of order $\ord$ is
a homomorphism from $\ca^*(V)$ to $\dz$ for any open $V \subset C$
(see, for example, proposition~\ref{lemma5} below). Then we have
$\ca^*/\ca_0^*\subseteq \dz_C$.

Let the sheaf $\ca^*/\ca_0^*|_U$ be a constant sheaf for an open set
$U$, which is equal to $d\dz$. (We suppose that this U is maximal.)
We have the following.
\begin{enumerate}
\item
\label{a6)}  If $H^0(C, \ca^*/\ca_0^*)=md\dz$ for some $m\ne 0$,
then $H^1(C, \ca^*/\ca_0^*)$ is a finite abelian group of order less
or equal to $m^{s-1}$ if $s> 1$, and is equal to 0 otherwise.
\item
\label{b6)}  If $H^0(C, \ca^*/\ca_0^*)=0$, then $rk (H^1(C,
\ca^*/\ca_0^*))\le {s-1}$ if $s> 1$, and  $H^1(C, \ca^*/\ca_0^*)=0$
otherwise.
\end{enumerate}
In both cases, $s$ is the number of critical points of
$\ca^*/\ca_0^*$, i.e. $s = \sharp (C \setminus U) $.
\end{prop}

\Proof. If $a\in \ca_x^*\cap (\ca_{j,x}\backslash \ca_{j+1,x})$,
where $x\in C$ is a point, and $b\in \ca_x^*$ is the inverse of $a$,
then $b\in\ca_x^*\cap (\ca_{-j,x}\backslash\ca_{-j+1,x})$. Then
$\ca_{j,x}=\ca_{0,x}a$. The relations
$\ca_{j,x}=\ca_{0,x}a$,
 $\ca_{-j,x}=\ca_{0,x}b$ and $ab=1$ extend to a neighbourhood $U$ of $x$.

Since $\ca_j/\ca_{j+1}$ is a torsion free sheaf, we obtain that if
$a\in \ca^*(U)$, then there exists a unique $j\in\dz$ such that
$a_x\in \ca_{j,x} \backslash \ca_{j+1,x}$, and the inverse $b$
satisfies $b_x\in \ca_{-j,x}\backslash \ca_{-j+1,x}$, and
$\ca_j|_U=\ca_0|_Ua$, $\ca_{-j}|_U=\ca_0|_Ub$. So we get in this
case an injection
$$
\ca^*/\ca_0^*\rightarrow \dz_C,\mbox{\quad} a\mapsto j=\ord (a)
$$
($\dz_C$ is a constant sheaf on $C$).
 This is an isomorphism iff $\ca_1,\ca_{-1}$ are mutually dual invertible $\ca_0$-modules.

By our assumptions we have that either $\ca^*=\ca_0^*$, or there is
a smallest positive integer $d$ such that there exists a point $x$
and $a\in\ca_x^*$ of order $d$. Then there is a largest open set $U$
where $\ca_d,\ca_{-d}$ are invertible mutually dual.

Then $\ca^*/\ca_0^*\subset d\dz$ and the cokernel is a sheaf with
support in $C\backslash U$. If $H^0(C, \ca^*/\ca_0^*)=0$, then at
least one stalk of the sheaf $d\dz / (\ca^*/\ca_0^*)$ in these
points is $d\dz$. If $H^0(C, \ca^*/\ca_0^*)=md\dz$, the stalks of
the sheaf $d\dz / (\ca^*/\ca_0^*)$ in these points are finite
groups, whose order is less or equal to $m$. Now, using the long
cohomology sequence of the short sequence
$$
0\rightarrow \ca^*/\ca_0^* \overset{\mu}{\rightarrow} d\dz \rightarrow coker(\mu ) \rightarrow 0
$$
we obtain the proof. (We use that  the first cohomology of a
constant sheaf on an irreducible space is equal to zero in Zariski
topology.)
\begin{flushright}
$\square$
\end{flushright}

\begin{prop}
\label{lemma5} Let $\xo_{\infty}$ be a ribbon with an irreducible
underlying curve $C$ over a field $k$. Assume that there exists a
point $x\in C$ such that $\ca_{1,x}\ca_{-1,x}=\ca_{0,x}$. Then the
function of order $\ord$ is compatible with the restriction
homomorphisms $\ca^*(U )  \to \ca^*(V) $ for open $V \subset U$,
and the function of order $\ord_U$ is a homomorphism from $\ca^*(U)$
to $\dz$ for any open $U$.
\end{prop}
\Proof. As it was shown in the proof of proposition \ref{prop6},
there is an invertible element $a\in \ca_{1,x}\backslash \ca_{2,x}$
such that $a^{-1}\in \ca_{-1,x}$. So, there exists an open $U\ni x$
such that $a\in \ca_1(U)$, $a^{-1}\in \ca_{-1}(U)$.

Now we need the following lemma.
\begin{lemma}
\label{lemma3} We consider a ribbon $(C, \ca )$, where $C$ is an
irreducible curve over a field $k$. Let the sheaf $\ca$ satisfies
the condition (*) (see definition~\ref{uslovie})  with the following
extra property: for any open $U_{\alpha}$ from (*) there exists an
invertible section $a\in \ca_1(U_{\alpha})\backslash
\ca_2(U_{\alpha})$ such that $a^{-1}\in \ca_{-1}(U_{\alpha})$.

Then the function of order $\ord$ is
compatible with the restriction homomorphisms $\ca^*(U )  \to
\ca^*(V) $ for open $V \subset U$,  and the function of order
$\ord_U$ is a homomorphism from $\ca^*(U)$ to $\dz$ for any open
$U$.
\end{lemma}

\Proof. The first assertion of the lemma follows from the second
one. Indeed, if $V\subset U$ are two open subsets and $b \in
\ca^*(U)$, $\ord_U(b)=k$, then $\ord_U(b^{-1})=-k$.
 We always have $\ord_V(b|_V)\ge \ord_U(b)$. If we suppose that $\ord_V(b|_V)> \ord_U(b)$,
 then $\ord_V((b|_V)^{-1})<-k=\ord_U(b^{-1})$. But $(b|_V)^{-1}=b^{-1}|_V$
 and $\ord_V(b^{-1}|_V)\ge \ord_U(b^{-1})=-k$, we have a contradiction.

Now we prove the second assertion of the lemma. At first, we prove
it for each $U_{\alpha}$. We  note that for any $b\in
\ca^*(U_{\alpha})$ and any $k\in \dz$
 we have $\ord_{U_{\alpha}}(ba^k)=\ord_{U_{\alpha}}(b)+k$,
 where $a$ is an invertible element from $\ca_1(U_{\alpha})\backslash \ca_2(U_{\alpha})$
such that $a^{-1} \in \ca_{-1}(U_{\alpha})$.
 Indeed, by definition of a ribbon,
 we always have $\ord_{U_{\alpha}}(bc)\ge \ord_{U_{\alpha}}(b)+\ord_{U_{\alpha}}(c)$
 for any $b,c \in \ca^*(U_{\alpha})$.
 Let $\ord_{U_{\alpha}}(ba)>\ord_{U_{\alpha}}(b)+1$. Then
$$
\ord_{U_{\alpha}}(b)=\ord_{U_{\alpha}}(baa^{-1})\ge
\ord_{U_{\alpha}}(ba)-1>\ord_{U_{\alpha}}(b) \mbox{,}
$$
we have a contradiction.

We note that $\ord_{U_{\alpha}}(bc)
=\ord_{U_{\alpha}}(b)+\ord_{U_{\alpha}}(c)$ if
$\ord_{U_{\alpha}}(c)=0$. Indeed, if $\ord_{U_{\alpha}}(bc)>
\ord_{U_{\alpha}}(b)$, then this would mean that $\bar{b}\bar{c}=0$,
where $\bar{b}\in \ca_{\ord (b)}(U_{\alpha})/\ca_{\ord
(b)+1}(U_{\alpha})$, $\bar{c}\in \co_C(U_{\alpha}) \ca_0(U_{\alpha})/ \ca_1(U_{\alpha})$.
But $\bar{c}, \bar{b}\neq 0$,
and $\ca_{\ord (b)}/\ca_{\ord (b)+1}$ is a torsion free sheaf by
definition, therefore we obtain a contradiction.

For any $b,c \in \ca^*(U_{\alpha})$ we have
$$
\ord_{U_{\alpha}}(bc)=\ord_{U_{\alpha}}(ba^{-\ord (b)}a^{\ord
(b)}c)=\ord_{U_{\alpha}}(ba^{-\ord (b)})+ \ord_{U_{\alpha}}(a^{\ord
(b)}c)= \qquad \qquad \qquad \qquad \qquad \qquad \qquad \qquad
$$
$$
\qquad \qquad \qquad \qquad \qquad \qquad \qquad \qquad \qquad
\qquad \qquad \qquad \qquad
 \ord_{U_{\alpha}}(b)+\ord_{U_{\alpha}}(c).
$$

The arguments from the beginning of the proof show that for any open $V\subset U_{\alpha}$ $\ord_V(a|_V)=1$,
 and $\ord_V((a|_V)^{-1})=-1$.
 Therefore, $\ord_V$ is also a homomorphism on $\ca^*(V)$.

Now let $U$ be an arbitrary open nonempty subset of $C$. Then
$U=\cup_{\alpha} (U\cap U_{\alpha})$, and $\ord_{U\cap U_{\alpha}}$
is a homomorphism for each $\alpha$. Let $b\in \ca^*(U)$,
$\ord_U(b)=k$. Assume that there exists $\beta$ such that
$\ord_{U\cap U_{\beta}}(b|_{U\cap U_{\beta}})=l>k$. Then for any
$\alpha $ we have $U\cap U_{\beta}\cap U_{\alpha}\neq \varnothing$
and $ \ord_{U\cap U_{\beta}\cap U_{\alpha}}(b|_{U\cap U_{\beta}\cap
U_{\alpha}})=l $ and therefore $\ord_{U\cap U_{\alpha}}(b|_{U\cap
U_{\alpha}})=l$. Since $\ca_{l}$ is a subsheaf of $\ca$, this would
mean that $b\in \ca_l(U)$, we have a contradiction.

So, we have for any $b,c\in \ca^*(U)$
$$
\ord_{U}(bc)=\ord_{U\cap U_{\alpha}}((bc)|_{U\cap U_{\alpha}})\ord_{U\cap U_{\alpha}}(b|_{U\cap U_{\alpha}})+\ord_{U\cap
U_{\alpha}}(c|_{U\cap U_{\alpha}})=\ord_U(b)+\ord_U(c) \mbox{.}
$$
The lemma is proved.
\begin{flushright}
$\square$
\end{flushright}

By lemma the function of order is a homomorphism on $U$ and on all open subsets of $U$.
Let $V\subset C$, $V\neq C$, $V\nsubseteq U$  be an open set. Since
$C$ is a reduced irreducible curve, $V$ must be affine and $V\cap U$
is also affine. Without loss of generality we can assume that $V\cap
U= D(f')$, where $V=\underline{\Spec }(B)$, $f'\in B$, $B=\co_C(V)$.
Let $f$ be a representative of $f'$ in $\ca_0(V)$. Clearly, it is
invertible in $\ca_0(V\cap U)$.   Let $b=a|_{V\cap U}$. We know that
$b$ is invertible, $\ord (b)=1$, $\ord (b^{-1})=-1$. Since the
sheaves $\ca_1/\ca_2$, $\ca_{-1}/\ca_0$ are coherent and $V\cap U$
is affine, there exists natural $n$ such that
$$
f^nb \mod  \ca_2(V\cap U)= \phi_{VD(f')}(\bar{b})\mbox{,} \qquad
f^nb^{-1} \mod  \ca_0(V\cap U)\phi_{VD(f')}(\overline{b^{-1}})\mbox{,}
$$
where $\phi_{VD(f')}: \ca (V)\rightarrow \ca (D(f'))$ is the
restriction homomorphism and $\bar{b}\in \ca_1/\ca_2(V)$,
$\overline{b^{-1}}\in \ca_{-1}/\ca_0(V)$, by \cite[lemma 5.3]{Ha}
and by proposition \ref{prop3}. Let $\tilde{b}, \widetilde{b^{-1}}$
be representatives of $\bar{b}, \overline{b^{-1}}$ in $\ca_1(V)$,
$\ca_{-1}(V)$ correspondingly. Then we have
$\ord_V(\tilde{b}\widetilde{b^{-1}})\ge 0$ and
$$\ord_V(\tilde{b}\widetilde{b^{-1}})\le \ord_{V\cap U}(\phi_{VD(f')}(\tilde{b}\widetilde{b^{-1}}))$$
by the properties of $\ca$. But
$\phi_{VD(f')}(\tilde{b}\widetilde{b^{-1}})= f^{2n} \mod \ca_1(V\cap
U)$, wherefrom $\ord_V(\tilde{b}\widetilde{b^{-1}})0$.

Note that for any $d\in \ca (V)$ we have
$\ord_{V}(dc)=\ord_{V}(d)+\ord_{V}(c)$ if $\ord_{V}(c)=0$. Indeed,
if $\ord_{V}(dc)> \ord_{V}(d)$, then this would mean that
$\bar{d}\bar{c}=0$, where $\bar{d}\in \ca_{\ord (d)}(V)/\ca_{\ord
(d)+1}(V)$, $\bar{c}\in \co_C(V)$. But the curve $C$ is reduced and
irreducible and $\ca_{\ord (d)}/\ca_{\ord (d)+1}$ is a torsion free
sheaf by definition, wherefrom we obtain a contradiction.

Now, repeating the arguments from the proof of lemma \ref{lemma3},
we obtain $\ord_V(d\tilde{b}^k)=\ord_V(d)+k$ for any integer $k$ and for any $d,c\in \ca^*(V)$
$$
\ord_{V}(dc)=\ord_{V}(d\widetilde{b^{-1}}^{\ord (d)}\tilde{b}^{\ord (d)}c)\ord_{V}(d\widetilde{b^{-1}}^{\ord (d)})+
\ord_{V}(\tilde{b}^{\ord (d)}c)=\ord_{V}(d)+\ord_{V}(c).
$$

At last, if $V=C$, then we can apply the arguments at the end of the
proof of lemma~\ref{lemma3}.\\
The proposition is proved.
\begin{flushright}
$\square$
\end{flushright}

\begin{corol}
\label{corol4} If there exists a point $P$ on an irreducible
curve $C$  such that $\ca_{1,P}\ca_{-1,P}=\ca_{0,P}$, then the
following properties are satisfied.
\begin{enumerate}
\item
\label{15)}  The embedding of sheaves
$\ca^*/\ca_0^*\overset{\ord}{\rightarrow}\dz_C$ is an isomorphism on
an open subset $U\subset C$. Besides, in the remaining points of
$C\backslash U$, the stalks $(\ca^*/\ca_0^*)_x$ are cyclic subgroups
$d_x\dz$ of $\dz$. If $H^0(C, \ca^*/\ca_0^*)=d\dz$ with $d>0$, then
all $d_x$ are divisors of $d$.
\item
\label{25)}  If $P$ is a smooth point on the curve $C$, then
$\dim_{k(P)} (\cm_P/\cm_P^2)=2$.
\end{enumerate}
\end{corol}
\Proof. The proof of assertion~\ref{15)} of this corollorary is
clear.

Now we prove assertion~\ref{25)}. From proposition~\ref{prop6} we
know that in our case $\ca_{i,P}=\ca_{0,P}a^i$ for all $i\ge 1$.
Since $P$ is a smooth point, we have $\cm_{C,P}=\co_{C,P}\bar{u}$
for some $\bar{u}\in \cm_{C,P}$. Let $u\in \ca_{0,P}$ be a
representative of $\bar{u}$. Then, clearly, $u,a$ generate the ideal
$\cm_P$ in the ring $\ca_{0,P}$ and are linearly independent in
$\cm_P/\cm_P^2$. So, we conclude that $\dim_{k(P)}
(\cm_P/\cm_P^2)=2$.
\begin{flushright}
$\square$
\end{flushright}

\begin{ex}
If a curve $C$ is not irreducible, then it is possible that the
function of order  is not a homomorphism from $\ca(U)^*$ to $\dz$
for open $U \subset C$.

For example, if we take an algebraizable ribbon from example \ref{ex1},
where $X$ is an affine plane and $C$ is a curve given by the equation $xy=0$,
 then the elements $x$ and $y$ will be invertible of order zero elements for any open $U \subset C$
 such that $U$ contains the point $(x=0, y =0)$.
 For, $(xy)$ is an invertible element from $\ca (U)$,
 and therefore $x^{-1} =y(xy)^{-1}$, $y^{-1}=x(xy)^{-1}$. But $\ord_U (xy)=1$, so, $\ord_U$ is
 not a homomorphism.
\end{ex}

\begin{ex}
Let $\xo_{\infty}$ be a ribbon from example \ref{ex1}, where $X$ is
assumed to be a smooth projective surface. Assume also that $(C\cdot
C)\neq 0$, and $C$ is an irreducible curve. We have that the
condition $\ca_{1,P}\ca_{-1,P}=\ca_{0,P}$ of corollary \ref{corol4}
is satisfied at each point $P \in C$. Therefore, by proposition
\ref{prop4} and corollary \ref{corol4} we have the following exact
sequence of sheaves on $C$:
$$
1 \lto \ca^*_0 \lto \ca^*  \lto \dz \lto 0  \mbox{,}
$$
and $H^0(C, \ca^*/\ca_0^*)=\dz$, $H^1(C, \ca^*/\ca_0^*)=0$. It gives
the following   exact sequence:
$$
0\rightarrow \dz \overset{\alpha}{\rightarrow} Pic (X_{\infty})\rightarrow Pic (\xo_{\infty})\rightarrow 0,
$$
where $\alpha (1)=\ca_1$. (The map $\alpha$ is an injective map,
because $\alpha(1)$ is not a torsion element in the group
$Pic(X_{\infty})$. Indeed, the image of $\alpha(1)$ in $Pic(C)$ has
degree equal to $- (C \cdot C) \ne 0$.) So, we obtain that $Pic
(\xo_{\infty})\simeq Pic (X_{\infty})/\langle \ca_1\rangle \simeq
Pic (X_{\infty})/\dz$.

For each $i$ we have the exact sequence
$$
0\rightarrow H^1(C, \frac{1+\ca_1}{1+ \ca_{i+1}})\rightarrow Pic
(X_i)\rightarrow Pic (C)\rightarrow 0
$$
and therefore we have the map
$$
Pic (X_i)\overset{deg}{\rightarrow} \dz \rightarrow 0, \mbox{\quad} \cl \mapsto deg (\cl |_C).
$$

By our assumptions we have $deg (\ca_1/\ca_{i+1})=d=-(C\cdot C)\ne 0$.
Therefore, we have the following exact diagrams for each $i$:
$$
\begin{array}{ccccccc}
0\rightarrow &Pic^0(X_i)&\rightarrow &Pic (X_i)& \rightarrow & \dz &\rightarrow 0\\
&\uparrow && \bigcup &&\bigcup & \\
&0 &\rightarrow & \langle \ca_1/\ca_{i+1}\rangle &\simeq &d\dz &\\
&&&\uparrow &&\uparrow &\\
&&&0 && 0&
\end{array}
$$
whence
$$
0\rightarrow Pic^0(X_i)\rightarrow Pic (X_i)/ \langle \ca_1/\ca_{i+1}\rangle \rightarrow \dz /d\dz \rightarrow 0.
$$

Projective system $(Pic(X_i), i \ge 0)$ satisfies the ML-condition
(as the first cohomology on the curve), therefore $(Pic^0(X_i), i
\ge 0)$ satisfies the ML-condition (as the kernels of the maps to
$\dz$). Passing to the projective limit we obtain the exact sequence
$$
\begin{array}{ccccccc}
0 \rightarrow & Pic^0(X_{\infty}) & \rightarrow & Pic(X_{\infty}) /
\langle \ca_1\rangle & \rightarrow & \dz /d\dz & \rightarrow 0 \\
&&&
\Vert \\
&&& Pic(\xo_{\infty})
\end{array}
$$

In particular, when $X=\dpp^2$, $C=\dpp^1\subset X$, we have $d=-1$ and
therefore $Pic^0(X_{\infty})\simeq Pic (\xo_{\infty})$.
On $Pic^0(X_{\infty})$ there exists a structure of a scheme (see, for example, \cite{Lip}),
so in this case there is a structure of a scheme also on $Pic (\xo_{\infty})$.
\end{ex}

\section{A generalized Krichever-Parshin map}
\label{KPmap}

Let $\xo_{\infty}=(C,\ca )$ be a ribbon over a field $k$.

\begin{defin}
\label{defin4} We say that a point $P\in C$ is a smooth point of the
ribbon $\xo_{\infty}$ if the following conditions are satisfied.
\begin{enumerate}
\item
\label{81)}  $P$ is a smooth point of $C$.
\item
\label{82)}  $\widehat{(\ca_i/\ca_{i+1})}_P\otimes
\widehat{(\ca_j/\ca_{j+1})}_P\rightarrow
\widehat{(\ca_{i+j}/\ca_{i+j+1})}_P$ is an isomorphism of
$\widehat{\co}_{C,P}$-modules, and this map is induced by the map
from the definition of ribbon: $\ca_i\cdot\ca_j\subset \ca_{i+j}$.
\end{enumerate}
\end{defin}

\begin{ex}
Let $\xo_{\infty}$ be a ribbon from example \ref{ex1}, where $P \in
C$ is a smooth point on the curve $C$ and the surface $X$. Then it
is clear that $P$ is a smooth point of the ribbon $\xo_{\infty}$.
\end{ex}

\begin{rem}
Immediately from definition follows that a ribbon with a smooth
point, whose topological space is irreducible, satisfies the conditions of proposition \ref{lemma5}. On the other
hand, all the ribbons from examples~\ref{ex2} and~\ref{ex2a} have
open neighborhoods, where they have smooth points  and are
algebraizable.
\end{rem}

\begin{prop}
\label{prop8}
Let $P$ be a smooth $k$-point of the ribbon $(C,\ca )$. Then
$$\tilde{\ca}_{0,P}\simeq \hat{\ca}_{0,P}\simeq k[[u,t]] \mbox{,}$$
where $t\tilde{\ca}_{0,P}=\tilde{\ca}_{1,P}$ and $\widehat{\co}_{C,P}\simeq k[[\tau (u)]]$, where
$\tau :\tilde{\ca}_{0,P}\rightarrow \widehat{\co}_{C,P}$ is a canonical map.
\end{prop}

\Proof. The isomorphism $\tilde{\ca}_{0,P}\simeq \hat{\ca}_{0,P}$
follows from corollaries~\ref{corol4} and~\ref{coroll}. Now we will
prove that $\widehat{\co}_{X_i,P}\simeq k[[u]][t]/t^i$ for some
$u,t$. The proof is by induction on $i$.

If $i=1$, then $\widehat{\co}_{X_1,P}=\widehat{\co}_{C,P}\simeq k[[u]]$ for some $u$.
Suppose we have proved the assertion for $(i-1)$. We have the exact triple:
$$
0\rightarrow \widehat{(\ca_{i-1}/\ca_i)}_P\rightarrow \widehat{\co}_{X_i,P}
\overset{\gamma}{\rightarrow}\widehat{\co}_{X_{i-1},P}\rightarrow 0.
$$
Let $\tilde{u}, \tilde{t}\in \widehat{\co}_{X_i,P}$ be elements with
$\gamma (\tilde{u})=u$, $\gamma (\tilde{t})=t$. From the definition
of a smooth point it follows that $\tilde{t}^{i-1}$ is a generator
of the $\widehat{\co}_{C,P}$-module $\widehat{(\ca_{i-1}/\ca_i)}_P$.
Therefore, $\widehat{\co}_{X_i,P}\simeq k[[u]][t]/t^i$.

Passing to the projective limit with respect to $i$ we are done.
\begin{flushright}
$\square$
\end{flushright}

\begin{defin}
Any elements $u,t$ from proposition \ref{prop8} are called
{\it formal local parameters} of the ribbon $(C,\ca )$ at the smooth point $P$.
\end{defin}

\begin{lemma}
\label{lemma4} Let $(C,\ca ,P, u, t)$ be a ribbon over a field $k$
with a smooth $k$-point $P$ and formal local parameters $u,t$. Then
$u\in \tilde{\ca}_{0,P}$ defines an effective Cartier divisor
$p_{u,i}$ on the scheme $X_i$ for any $i$ such that
$\theta_i^*p_{u,i}=p_{u,i-1}$ and $p_{u,1}=P$, where
$$
\theta_i: X_{i-1}\rightarrow X_i
$$
is a canonical map.
\end{lemma}

\Proof. We know by proposition \ref{prop8} that
$\widehat{\co}_{X_i,P}\simeq k[[u]][t]/t^i$, because $P$ is a smooth
point of the ribbon $(C,\ca )$. By $\tilde{p}_{u,i}=u\cdot
k[[u]][t]/t^i$ we denote the ideal in $\widehat{\co}_{X_i,P}$. Let
$p_{u,i}':=\tilde{p}_{u,i}\cap \co_{X_i,P}$ be the ideal in
$\co_{X_i,P}$.

We have for some $j>0$ that $\cm_P^j\cdot k[[u]][t]/t^i\subset
\tilde{p}_{u,i}$, where $\cm_P$ is the maximal ideal of
$\co_{X_i,P}$. Therefore, $\cm_P^j\co_{X_i,P}\subset p_{u,i}'$. Let
$\tilde{u} \in \co_{X_i,P}$ be an element such that $\beta
(\tilde{u})$ coincides with $\hat{\beta}(u)$, where $\beta,
\hat{\beta}$ are the following natural maps
$$
\begin{array}{rclcc}
\beta & \ : & \co_{X_i,P} & \longrightarrow & \co_{X_i,P}/\cm_P^j \\
      &     &             &                 &  \| \\
\hat{\beta} & : & \widehat{\co}_{X_i,P} & \lto &
\widehat{\co}_{X_i,P}/ \widehat{\cm}_{P}^j
\end{array}
$$
Then $\tilde{u}\in p_{u,i}'$ and $\tilde{u}\cdot
\widehat{\co}_{X_i,P}=\tilde{p}_{u,i}$. Therefore, $p_{u,i}'\cdot
\widehat{\co}_{X_i,P}=\tilde{p}_{u,i}$, and
$p_{u,i}'=\tilde{u}\co_{X_i,P}$ defines the effective Cartier
divisor $p_{u,i}$ in some affine open neigbourhood of the point
$P\in X_i$ (and on $X_i$). By construction,
$\theta^*_ip_{u,i}=p_{u,i-1}$.
\begin{flushright}
$\square$
\end{flushright}

\begin{rem}
By construction of the ideal $p_{u,i}'$ (or divisor $p_{u,i}$) we
obtain that it is uniquely defined by the properties $p_{u,i}'\cdot
\widehat{\co}_{X_i,P}=u\cdot \widehat{\co}_{X_i,P}$,
$\theta^*_ip_{u,i}=p_{u,i-1}$, and $p_{u,1}'=\cm_P\subset
\co_{C,P}$.
\end{rem}

\begin{defin} \label{tfsh}
Let $\xo_{\infty}=(C,\ca )$ be a ribbon over a scheme $S$. We say
that $\cn$ is a torsion free sheaf of rank $r$ on $\xo_{\infty}$ if
  $\cn$ is a sheaf of $\ca$-modules on $C$ with a descending
filtration $(\cn_i)_{i\in \sdz}$ of $\cn$ by $\ca_0$-submodules
which satisfies the following axioms.
\begin{enumerate}
\item
 $\cn_i\ca_j\subseteq \cn_{i+j}$ for any $i,j$.
\item \label{ittf}
For each $i$ the sheaf $\cn_i/\cn_{i+1}$ is a coherent sheaf on $C$,
flat over $S$, and for any $s \in S$  the sheaf
$\cn_i/\cn_{i+1}|_{C_S}$ has no coherent subsheaf with finite
support, and is isomorphic to $\co_{C_S}^{\oplus r}$ on a dense open
set.
\item
 $\cn =\liminj\limits_i \cn_i$ and
$\cn_i=\limproj\limits_{j>0}\cn_i/\cn_{i+j}$ for each $i$.
\end{enumerate}
\end{defin}

\begin{rem}
It follows from  assertion~\ref{ittf} of  definition~\ref{tfsh} that
if $C_s$ (for $s \in S$) is an irreducible curve, then the sheaf
$\cn_i /\cn_{i+1} \mid_{C_s}$ is a rank~$r$ torsion free sheaf on
$C_s$ for any $i \in\dz$.
\end{rem}

\begin{rem} \label{remtf}
If the sheaf $\ca$ of a ribbon $\xo_{\infty}$ satisfies the
condition (*) from definition~\ref{uslovie}, then any torsion free
sheaf $\cn$ of rank $r$ on $\xo_{\infty}$ is coherent. The proof of
this fact is the same as in lemma \ref{lemma2}.
\end{rem}

On the other hand, if $\ca$ is only coherent, then there exists a
torsion free sheaf that is not coherent, as it follows from the
example below.
\begin{ex}
\label{ex4} Consider the ribbon $(C, \ca =\co_C ((t))^{Q} )$ from
example \ref{ex2a}. The sheaf $\ca$ is coherent, but not weakly
Noetherian. The sheaf $\cn :=\co_C ((t))$ with obvious filtration is
a torsion free sheaf of rank 1 on $\xo_{\infty}$. But the stalk
$\cn_Q$ can not be finitely generated: for any finite number of
sections $g_1, \ldots ,g_k\in \cn (V)$, $Q\in V$ there are infinite
number of elements $t^{l}$, $l\ll 0$ that can not be generated by
$g_i$. So, $\cn$ is not of finite type and therefore is not
coherent.
\end{ex}

\begin{defin}
Let $\xo_{\infty}=(C,\ca )$ be a ribbon over a field $k$. We say
that a point $P\in C$ is a smooth point of a torsion free sheaf
$\cn$ on $\xo_{\infty}$ if the following conditions are satisfied.
\begin{enumerate}
\item
 $P$ is a smooth point of $\xo_{\infty}$.
\item
 $\widehat{(\cn_i/\cn_{i+1})}_P\otimes_{\hat{\co}_{C,P}}
\widehat{(\ca_j/\ca_{j+1})}_P\rightarrow
\widehat{(\cn_{i+j}/\cn_{i+j+1})}_P$ is an isomorphism of
$\widehat{\co}_{C,P}$-modules, and this map is induced by the map
from the definition of $\cn$: $\cn_i\cdot\ca_j\subset \cn_{i+j}$.
\end{enumerate}
\end{defin}

Similarly to the proposition~\ref{prop8} we have the following
proposition.
\begin{prop}
Let $P$ be a smooth point of a torsion free sheaf $\cn$ of rank $r$
on a ribbon $\xo_{\infty}$ over a field $k$. Then
$$\widetilde{\cn}_{0,P}\simeq \widetilde{\ca}_{0,P}^{\oplus r}
\mbox{,}$$
where $\widetilde{\cn}_{0,P}=\limproj\limits_{j\ge
0}\widehat{(\cn_0/\cn_j)}_P$.
\end{prop}

\Proof. By induction on $j$ and using the exact sequence
$$
0\rightarrow \widehat{(\cn_{j-1}/\cn_j)}_P\rightarrow \widehat{(\cn_0/\cn_j)}_P\rightarrow
\widehat{(\cn_0/\cn_{j-1})}_P\rightarrow 0
$$
we prove that $\widehat{(\cn_0/\cn_j)}_P\simeq
\widehat{(\ca_0/\ca_j)}_P^{\oplus r}$. Then we pass to the
projective limit.
\begin{flushright}
$\square$
\end{flushright}

\begin{ex}
Let $\xo_{\infty}= (C, \ca)$ be a ribbon such that $C$ is
irreducible. We suppose that the function of order  is a
homomorphism (see, for example, proposition~\ref{lemma5} above.)
Then every element $\cn \in Pic (\xo_{\infty})$  gives a torsion
free sheaf of rank $1$ on $\xo_{\infty}$ after the fixing of a
filtration on $\cn$. (All the possible filtrations form a
$\dz$-torsor.) Indeed, the function of order $\ord$ gives a
homomorphism:
$$
\gamma \quad : \quad H^1(C, \ca^*) \lto H^1(C, \dz) \mbox{.}
$$
And the obstruction to find a filtration on $\cn$ is $\gamma(\cn)
\ne 0$. But any local $\dz$-system on irreducible space $C$ is
trivial in Zariski topology, i.e. $H^1(C, \dz) =0$. Thus, we have a
filtration on $\cn$.

\end{ex}

\begin{ex}
 Let $\xo_{\infty}$ be a ribbon from example~\ref{ex1}. Let $E$ be a locally free sheaf of rank $r$
 on the surface $X$. Then
$$\eo_C:=\liminj\limits_i\limproj\limits_jE(iC)/E(jC)$$
is a torsion free sheaf of rank $r$ on $\xo_{\infty}$. Any point
$P\in C\subset X$ that is smooth on $C$ and on $X$ will be a smooth
point on $\eo_C$.
\end{ex}

\begin{rem}
Similarly to definition~\ref{defrin} we have two
$\ca_{0,P}$-modules: $\widetilde{\cn}_{0,P}$ and
$\widehat{\cn}_{0,P}$, where the latter is  the $\cm_P$-adic
completion of the module $\cn_{0,P}$. Using similar arguments as in
the the proof of proposition \ref{prop5}, we obtain that if
$\dim_k\cn_{0,P}/\cm_P\cn_{0,P} < \infty$ , then the natural
homomorphism of $\ca_{0,P}$-modules
$$
\widehat{\cn}_{0,P}\overset{\alpha}{\rightarrow} \widetilde{\cn}_{0,P}
$$
is surjective.

If $P$ is a smooth point of the torsion free sheaf $\cn$ of
rank~$r$, then $\widehat{\ca}_{0,P}\simeq \widetilde{\ca}_{0,P}$,
$\dim_{k( P)}\cn_{0,P}/\cm_P\cn_{0,P}=r$ and therefore $\alpha$ is
an isomorphism of $\widetilde{\ca}_{0,P}$-modules.
\end{rem}

Further we will work with geometric data consisting of a ribbon, a
torsion free sheaf on it, formal local parameters at a point $P$ on
the ribbon and a formal trivialization of the sheaf at $P$.

\begin{defin} \label{defin13}
Let $(\xo_{\infty}, \cn )$, $(\xo'_{\infty}, \cn' )$ be two ribbons
over a field $k$ with two torsion free sheaves of rank $r$ on them.
We say that $(\xo_{\infty}, \cn )$ is isomorphic to $(\xo'_{\infty},
\cn' )$ if there is an isomorphism
$$\varphi \; : \; \xo_{\infty} \longrightarrow
\xo'_{\infty}$$ of ribbons (see definition~\ref{defisom}) and an
isomorphism
 $$\psi \; : \; \cn'\rightarrow \varphi_{*}(\cn )$$ of graded $\ca'$-modules,
 i.e. $\psi (\cn'_i)=\varphi_{*}(\cn_i)$ and $\psi (ln)=\varphi^{\sharp}(l)\psi (n)$
 for any sections $n\in \cn' (U)$, $l\in \ca' (U)$ over open $U\subset C'$.
\end{defin}
\begin{defin}
\label{defin9} We will consider the following geometric data $(C,\ca
, \cn ,P, u,t, e_P)$, where
\begin{itemize}
\item
$(C,\ca )$ is a ribbon over a field $k$,
\item
$\cn$ is a torsion free sheaf of rank $r$ on $(C, \ca )$,
\item
$P\in C$ is a  smooth $k$-point of the sheaf $\cn$,
\item
$u,t$ are formal local parameters of the ribbon at $P$,
\item
$e_P: \widetilde{\cn}_{0,P} \rightarrow \widetilde{\ca}_{0,P}^{\oplus r}\simeq k[[u,t]]^{\oplus r}$
is an isomorphism of $\widetilde{\ca}_{0,P}$-modules.
\end{itemize}

We say that $(C,\ca , \cn ,P, u,t, e_P)$ is isomorphic to $(C',\ca'
, \cn' ,P', u',t', e'_P)$ if there is an isomorphism (see
definition~\ref{defin13})
$$
(\varphi ,\psi ) \; : \; (C,\ca ,\cn ) \longrightarrow (C',\ca'
,\cn' )
$$
such that $\varphi (P)=P'$, $\varphi_P^{\sharp}(t')=t$, $\varphi_P^{\sharp}(u')=u$,
where $\varphi_P^{\sharp}:\widetilde{\ca'}_{0,P}\rightarrow \widetilde{\ca}_{0,P}$
is an isomorphism of local rings induced by $\varphi^{\sharp}$, and the diagram
$$
\begin{array}{ccc}
\widetilde{\cn'}_{0,P}&\overset{\psi_P}{\rightarrow}&\widetilde{\cn}_{0,P}\\
\downarrow \lefteqn{e'_P}&&\downarrow \lefteqn{e_P}\\
\widetilde{\ca'}_{0,P}^{\oplus r}&\overset{\varphi_P^{\sharp}}{\rightarrow}&\widetilde{\ca}_{0,P}^{\oplus r}
\end{array}
$$
is commutative, where the isomorphism $\psi_P$ is induced by $\psi$.
\end{defin}


\begin{defin}
Let $K=k((u))((t))$ be a two-dimensional local field. We define the following $k$-subspaces of $K$:
$$
\co (n)=t^nk((u))[[t]]
$$
for any $n$. For any $k$-subspace $W\subset K^{\oplus r}$ and any
$j>i \in \dz$ we define
$$
W(i,j)=\frac{W\cap \co (i)^{\oplus r}}{W\cap \co (j)^{\oplus r}}
 \mbox{.}
$$
\end{defin}
We have the natural isomorphism $\co(i)/ \co(j) = k((u))^{\oplus
(j-i)}$, therefore $W(i,j)$ is a $k$-subspace of the space
$k((u))^{\oplus r(j-i)}$. We note that the last space has a natural locally linearly compact topology.
Also the spaces $K$ and $K^{\oplus r}$ have a natural topology of inductive-projective limits.

\begin{defin}
\label{defin11} Let $W$ be a closed $k$-subspace of $K^{\oplus
r}=k((u))((t))^{\oplus r}$, let $A$ be a closed $k$-subalgebra of
$K=k((u))((t))$. (We can consider $K^{\oplus r}$ as a $K$-module, so
the product $A\cdot W\subset K^{\oplus r}$ is defined.)

We suppose that $A\cdot W\subset W$, and $A(i,i+1)\subset k((u))$
is a discrete subspace with quotient being linearly compact space,
$W(i,i+1)\subset k((u))^{\oplus r}$
is a discrete subspace with quotient being linearly compact
for any $i\in\dz$. Then we call the pair of
$k$-subsapces $(A,W)\subset K\oplus K ^{\oplus r}$ as a {\it Schur
pair}.
\end{defin}

(It is easy to see that the property to be closed subspaces for $W$ and $A$ is equivalent to the following equalities:
$
W =\liminj\limits_i\limproj\limits_{j\ge i} W (i,j)$
and
$A =\liminj\limits_i\limproj\limits_{j\ge i} A (i,j)  \mbox{.}
$)

\begin{rem}
By induction on $j-i > 0$ we have that if $W(i,i+1)$, for any $i \in \dz$, are
subspaces from definition \ref{defin11},
then $W(i,j)$ is a discrete subspace in
$k((u))^{\oplus r(j-i)}$ with quotient being linearly compact
for any $j>i$. Similarly, $A(i,j)$
is a discrete subspace in $k((u))^{\oplus (j-i)}$ with quotient being linearly compact for any $j>i$.
\end{rem}

\begin{rem}
Clearly, the subspaces $W(i,i+1)$ (resp. $A(i,i+1)$) from
definition~\ref{defin11} satisfy the Fredholm condition with respect
to $k[[u]]^{\oplus r}$ (resp. to $k[[u]]$, see introduction), as the
analogous subspaces in the construction of the Krichever map in
\cite{Pa}, \cite{Os}, see also \cite{ZhO}.
\end{rem}

\begin{theo}
\label{th1} The Schur pairs $(A,W)$ from definition \ref{defin11}
are in one-to-one correspondence with  data $(C,\ca , \cn ,P, u,t,
e_P)$ from definition \ref{defin9} up to an isomorphism, where we
additionally assume that $C$ is a projective irreducible curve.
\end{theo}
\begin{corol}
$k$-subalgebras $A$ from definition \ref{defin11} are in one-to-one
correspondence with data $(C,\ca ,P, u,t)$ up to an isomorphism,
where $C$ has to be a projective irreducible curve.
\end{corol}

\Proof. The corollary follows from the theorem, if we take
$\cn=\ca$, $W=A$, and $e_P=1$.

Now we'll prove the theorem.  We have the following diagram of maps
for any coherent sheaf $M$ on the scheme $X_i$ and for any $i\ge 0$.
$$
\begin{array}{cccc}
\Gamma (X_i\backslash P, M)\overset{\alpha_M}{\rightarrow}&\Gamma
(\Spec  \co_{X_i,P}\backslash P, M)&\overset{\beta_M
}{\rightarrow}&\Gamma
(\Spec  \widehat{\co}_{X_i,P}\backslash P, M)\\
&\parallel &&\parallel \\
&M\otimes_{\co_{X_i}}\co_{X_i,\eta_i}&&M\otimes_{\co_{X_i}}
\co_{X_i,\eta_i}\otimes_{\co_{X_i}}\widehat{\co}_{X_i,P}
\end{array}
$$
where $\eta_i$ is the generic point of the scheme $X_i$.

Now let $M=\cn_k/\cn_{k+i+1}$ for some $k$. Then, by
statement~\ref{i1} of proposition \ref{prop3} and statement~\ref{2)}
of proposition~\ref{prop1}, $M$ is coherent sheaf on the scheme
$X_i$. By induction on $i$, we show  that the map $\alpha_M$ is an
embedding. It is true for $i=0$, because $\cn_k/\cn_{k+1}$ is a
torsion free sheaf on $C$. We have the following commutative diagram
for arbitrary $i \ge 1$
$$
\begin{array}{ccccccc}
0\rightarrow & \Gamma (X_i\backslash P, \cn_{i+k}/\cn_{i+k+1}) & \rightarrow &
\Gamma (X_i\backslash P, \cn_k/\cn_{i+k+1}) & \rightarrow & \Gamma (X_i\backslash P, \cn_k/\cn_{k+i}) &
\rightarrow 0\\
&&&&&&\\
& \downarrow \lefteqn{\alpha_{\cn_{i+k}/\cn_{i+k+1}}} &&
 \downarrow \lefteqn{\alpha_{\cn_{k}/\cn_{i+k+1}}} && \downarrow \lefteqn{\alpha_{\cn_{k}/\cn_{i+k}}} & \\
&&&&&&\\
0 \rightarrow & (\cn_{i+k}/\cn_{i+k+1})_{\eta_i} &\rightarrow &
(\cn_k/\cn_{i+k+1})_{\eta_i} &\rightarrow &
(\cn_k/\cn_{k+i})_{\eta_i}
&\rightarrow 0,\\
&&&&&\parallel & \\
&&&&&(\cn_k/\cn_{k+i})_{\eta_{i-1}}&
\end{array}
$$
since $\cn_{i+k}/\cn_{i+k+1}$ is a coherent $\co_{X_i}$-module and
the $\co_{X_i}$-module structure on $\cn_k/\cn_{k+i}$ is the same as
the $\co_{X_{i-1}}$-module structure. Therefore, by induction
hypothesis, the left and right vertical arrows are  embeddings.
Hence, the middle arrow is also an embedding.

The map $\beta_{\cn_k/\cn_{i+k+1}}$ is an embedding for the sheaf
$\cn_k/\cn_{i+k+1}$. Therefore, the map \linebreak
$\beta_{\cn_k/\cn_{i+k+1}} \circ \alpha_{\cn_k/\cn_{i+k+1}}$ is an
embedding for the sheaf $\cn_k/\cn_{i+k+1}$.

Now  we have for $k=0$
$$
\ca_0/\ca_{i+1}\otimes_{\co_{X_i}}
\co_{X_i, \eta_i}\otimes_{\co_{X_i}}\widehat{\co}_{X_i,P}\simeq k((u))[t]/t^{i+1},
$$
because we have fixed the formal local parameters $u,t$ of our ribbon at $P$.

We have for $k>0$
$$
\ca_k/\ca_{k+i+1}\otimes_{\co_{X_i}}\co_{X_i, \eta_i}\otimes_{\co_{X_i}}
\widehat{\co}_{X_i,P}\simeq t^k\cdot k((u))[t]/t^{i+1}
$$
as an ideal in $\ca_0/\ca_{i+1+k}\otimes_{\co_{X_i}}\co_{X_i, \eta_i}
\otimes_{\co_{X_i}}\widehat{\co}_{X_i,P}\simeq k((u))[t]/t^{i+k+1}$.

By definition of a smooth point on a ribbon, we have  the natural
pairing for $k<0$
$$
\widehat{(\ca_k/\ca_{k+i+1})}_P\otimes
\widehat{(\ca_{-k}/\ca_{-k+i+1})}_P\rightarrow
\widehat{(\ca_0/\ca_{i+1})}_P \mbox{,}
$$
and the element $t^{-k} \mod \widetilde{\ca}_{-k+i,P}\in
\widehat{(\ca_{-k}/\ca_{-k+i+1})}_P$. Then, by induction on $i \ge
0$, we obtain that
$$
\widehat{(\ca_k/\ca_{k+i+1})}_P \overset{\times t^{-k}}{\longrightarrow}
\widehat{(\ca_0/\ca_{i+1})}_P\simeq k((u))[t]/t^{i+1}
$$
is an isomorphism. Therefore, we have
$$
\liminj\limits_k\limproj\limits_{i \ge 0}
\widehat{(\ca_k/\ca_{k+i+1})}_P\simeq k((u))((t)).
$$
Besides, $\ca =\liminj\limits_k\limproj\limits_{i \ge 0}
\ca_k/\ca_{k+i+1}$. Therefore, the ring
$$
A:=\liminj\limits_k\limproj\limits_{i \ge 0} \Gamma (X_i\backslash
P, \ca_k/\ca_{k+i+1}) \subset k((u))((t))
$$
is a $k$-subalgebra that satisfies the conditions of the theorem.

Analogously, using the trivialization $e_P$ and formal local
parameters $u,t$, we obtain the isomorphism
$$
\liminj\limits_k\limproj\limits_{i \ge 0}
\widehat{(\cn_k/\cn_{k+i+1})}_P\simeq k((u))((t))^{\oplus r}
$$
and the subspace
$$
W:=\liminj\limits_k\limproj\limits_{i \ge 0} \Gamma (X_i\backslash
P, \cn_k/\cn_{k+i+1}) \subset k((u))((t))^{\oplus r}
$$
is a $k$-subspace that satisfies the conditions of the theorem.

Thus, starting from the geometric data $(C,\ca , \cn ,P, u,t, e_P)$
from definition \ref{defin9}, we have constructed a Schur pair
$(A,W)$ from definition \ref{defin11}.

Now we are going to construct a geometric data starting from a Schur
pair. At first, we note that
$$
\Gamma (\Spec  \co_{X_i,P}\backslash P, \ca_k/\ca_{k+i+1})=\liminj\limits_{n \ge 0}
\Gamma (X_i, \ca_k/\ca_{k+i+1}(np_{u,i})),
$$
where $p_{u,i}$ is the effective Cartier divisor on $X_i$ which was
constructed in lemma~\ref{lemma4} above.

We consider the $k$-subspaces for $j > i \in \dz$
$$A(i,j)\subset
k((u))^{\oplus (j-i)} \quad \mbox{and}$$
$$U_n(i,j)=u^{-n}\cdot
k[[u]]^{\oplus (j-i)}\subset k((u))^{\oplus (j-i)} \mbox{.}$$ If
$i=0$, the space $\bigoplus_{n \ge 0}(U_n(0,j)\cap A(0,j))$ is a
graded ring. We put
$$X_{j-1}=\Proj (\bigoplus_{n \ge 0}(U_n(0,j)\cap A(0,j))).$$
The image of the embedding of $\bigoplus_{n \ge 0}(U_{n-1}(0,1)\cap
A(0,1))$ in $\bigoplus_{n \ge 0}(U_{n}(0,1)\cap A(0,1))$ is a
homogeneous ideal that determines a point $P\in X_0$.

If $j >i \in \dz$, then $\bigoplus_{n\ge 0}(U_n(i,j)\cap A(i,j))$ is
a graded module over the graded ring $\bigoplus_{n \ge
0}(U_n(0,j-i)\cap A(0,j-i))$. Then we define
$$\ca (i,j)=\widetilde{
\bigoplus_{n\ge 0}(U_n(i,j)\cap A(i,j))} \mbox{,}$$ i.e. it is a
coherent sheaf on $X_{(j-i)}$ which is associated with the
corresponding graded module. Since  there is no zero divisors in the
field $k((u))$, the sheaf $\ca (i,i+1)$ is a torsion free sheaf on
$C$ for any $i$.

For all $j > i \in \dz$ we have surjective morphisms $\ca (i,j+1)
\rightarrow \ca (i,j)$ and injective morphisms $\ca (i,j)
\rightarrow \ca (i-1,j)$. Also, from definitions, we have  maps for
all $i < j$, $k < l$
\begin{equation} \label{mulmap}
A(i,j)\otimes_k A(k,l)\longrightarrow A(i+k, \min (j+k,i+l))
\mbox{,}
\end{equation}
which are also  well-defined maps, if we pass to projective limits
with respect to $j$ and $l$.

So, we define sheaves $\ca$, $\ca_i$, $i \in \dz$ by
$$
\ca :=\liminj\limits_i\limproj\limits_{j\ge i}\ca (i,j), \mbox{\quad} \ca_i=\limproj\limits_{j\ge i}\ca (i,j).
$$
The map given by formula~(\ref{mulmap}) defines the multiplication
$\ca_i\cdot \ca_j\subset \ca_{i+j}$. Besides,
$\widetilde{\ca}_{0,P}=k[[u,t]]$, and therefore $u,t$ are the formal
local parameters of the ribbon $(X_0, \ca )$ at the point $P$.

Analogously, we define sheaves of modules $\cn$, $\cn_i$, $i \in
\dz$ by
$$
\cn :=\liminj\limits_i\limproj\limits_{j\ge i}\cn (i,j)\mbox{,}
\qquad \cn_i=\limproj\limits_{j\ge i}\cn (i,j) \mbox{,}
$$
where $\cn (i,j)=\widetilde{N} $, $N= \bigoplus_{n \ge
0}((u^{-n}\cdot k[[u]]^{\oplus r(j-i)})\cap W(i,j))$, i.e.
$\cn(i,j)$ is a coherent sheaf of $\co_{X_{(j-i-1)}}$-modules, which
is associated with the corresponding graded module, for all $j > i$
. By construction, we have a natural isomorphism
$$e_P \;: \; \widetilde{\cn}_{0,P}\rightarrow
k[[u,t]]^{\oplus r} \mbox{.}$$

The map $(A,W) \mapsto (X_0, \ca , \cn ,P, u,t, e_P)$ just
constructed is the inverse to the map which was constructed in the
first part of  proof of this theorem, because the sheaf
$\cn_i/\cn_{j}\simeq \widetilde{\Gamma_*(\cn_i/\cn_{j+1})}$ for all
$j > i \in \dz$, where the graded module
$$\Gamma_*(\cn_i/\cn_{j})=\bigoplus_{n \ge 0} \Gamma (X_{j-i-1}, \cn_i/\cn_{j}(np_{u,j-i-1})) \mbox{,}$$
defines a coherent sheaf  on the scheme
$$X_{j-i-1} = \Proj
(\bigoplus_{n\ge 0}\Gamma (X_{j-i-1},
\co_{X_{j-i-1}}(np_{u,j-i-1})))  \mbox{,} $$
 since
$\co_{X_{j-i-1}}(p_{u,j-i-1})$ is an ample sheaf on $X_{j-i-1}$. The
latter follows from the  lemma below.
\begin{lemma}
For any $i>0$ the sheaf $\co_{X_i}(p_{u,i})$ is an ample sheaf on $X_i$.
\end{lemma}

\Proof . $X_i$ is a proper scheme (as $X_0$ is a projective curve).
So, by \cite[ch.III, prop.5.3]{Ha} it is enough to prove that for
any $l>0$, for any coherent sheaf $\cf$ on $X_i$ there exists
$n_0>0$ such that for any $n>n_0$ $H^l(X_i, \cf \otimes
\co_{X_i}(np_{u,i}))=0$.

We use induction on $i$. If $i=1$, then $p_{u,1}$ is the point $P$
on the projective curve $C=X_0$, i.e. it is an ample divisor. If
$i>1$, we consider the exact sequence of $\co_{X_i}$-modules
\begin{multline*}
0\rightarrow \cf\otimes_{\co_{X_i}}\ca_{i-1}/\ca_i \otimes_{\co_{X_i}}\co_{X_i}(np_{u,i}) \rightarrow
\cf \otimes_{\co_{X_i}}\co_{X_i}(np_{u,i})\rightarrow \\
\cf\otimes_{\co_{X_i}} (\ca_0/\ca_{i-1})\otimes_{\co_{X_i}}\co_{X_i}(np_{u,i})  \rightarrow 0.
\end{multline*}

The $\co_{X_i}$-module structure of modules
$\cf\otimes_{\co_{X_i}}\ca_{i-1}/\ca_i$ and $\cf\otimes_{\co_{X_i}}
(\ca_0/\ca_{i-1})$ coincides with the $\co_{X_{i-1}}$-module
structure. Therefore, their cohomology on $X_i$ coincide with
cohomology on $X_{i-1}$. So, from the long exact cohomology sequence
and induction hypothesis we get for all $n>n_0$ and all $l>0$
$H^l(X_i, \cf \otimes \co_{X_i}(np_{u,i}))=0$.
\begin{flushright}
$\square$
\end{flushright}
The theorem is proved.
\begin{flushright}
$\square$
\end{flushright}

\begin{rem}
The constructions of subspaces and geometric data given in the
theorem are generalizations of the Krichever map constructed in the
works \cite{Pa}, \cite{Os}.
 If a geometric data is taken on a ribbon that comes from a surface
 and a reduced effective Cartier divisor on it, as in example \ref{ex1},
 then these constructions coincide.
\end{rem}

\noindent H. Kurke,  Humboldt University of Berlin, department of
mathematics, faculty of mathematics and natural sciences II, Unter
den Linden 6, D-10099, Berlin, Germany \\ \noindent\ e-mail:
$kurke@mathematik.hu-berlin.de$

\vspace{0.5cm}

\noindent D. Osipov,  Steklov Mathematical Institute, algebra
department, Gubkina str. 8, Moscow, 119991, Russia \\ \noindent
e-mail:
 ${d}_{-} osipov@mi.ras.ru$

\vspace{0.5cm}

\noindent A. Zheglov,  Lomonosov Moscow State  University, faculty
of mechanics and mathematics, department of differential geometry
and applications, Leninskie gory, GSP, Moscow, \nopagebreak 119899,
Russia
\\ \noindent e-mail
 $azheglov@mech.math.msu.su$

\end{document}